\documentclass[english, 12pt,leqno]{amsart}

\usepackage{amsmath}
\usepackage{amsfonts}
\usepackage{amssymb}
\newtheorem{Def}{D\'efinition}
\newtheorem{Thm}[Def]{Theorem}
\newtheorem{Prop}[Def]{Proposition}
\newtheorem{Lem}[Def]{Lemma}
\newtheorem{Cla}[Def]{Claim}
\newtheorem{Cor}[Def]{Corollary}

\newtheorem{Ahah}[Def]{Homotopy Annulus Theorem}

\title[On malnormal peripheral subgroups]
{On malnormal peripheral subgroups 
\\
in fundamental groups of $3$-manifolds}

\author{Pierre de la Harpe and Claude Weber}

\date{April 10,  2011}

\address{
Section de math\'ematiques, 
Universit\'e de Gen\`eve, 
C.P.~64, 
CH--1211 Gen\`eve 4. 
\hskip.4cm
Pierre.delaHarpe@unige.ch
\hskip.2cm \& \hskip.2cm
Claude.Weber@unige.ch
\newline
}



\begin{document}

\begin{abstract}
Let $K$ be a non-trivial knot in the $3$-sphere, 
$E_K$ its exterior,
$G_K = \pi_1(E_K)$ its group, 
and $P_K = \pi_1(\partial E_K) \subset G_K$ 
its peripheral subgroup.
We show that  $P_K$ is malnormal in $G_K$, 
namely that $gP_Kg^{-1} \cap P_K = \{e\}$ for any $g \in G_K$ with $g \notin P_K$,
unless $K$ is in one of the following three classes:
torus knots, cable knots, and composite knots;
these are exactly the classes for which
there exist annuli in $E_K$ attached to $T_K$
which are not boundary parallel
(Theorem \ref{laSimon} and Corollary \ref{C}).
More generally, we characterise malnormal peripheral subgroups
in the fundamental group of a compact orientable irreducible $3$-manifold 
with boundary a non-empty union of tori (Theorem \ref{Thm3manif}).
Proofs are written with non-expert readers in mind.
Half of our paper (Sections \ref{section7} to \ref{sectionDigression})
is a reminder of some three-manifold topology 
as it flourished before the Thurston revolution.
\par
In a companion paper \cite{HaWeOs}, 
we collect general facts on malnormal subgroups and Frobenius groups,
and we review a number of examples.
\begin{center}
Table of sections.
\end{center}
\par\noindent
1. Statement of the results.
\par\noindent
2. On annuli embedded in $3$-manifolds
with torus boundaries.
\par\noindent
3. Examples of non boundary-parallel annuli in knot exteriors.
\par\noindent
4. Proof of Theorem \ref{laSimon}.
\par\noindent
5. Proof of Corollary \ref{C} and of Theorem \ref{Thm3manif}.
\par\noindent
6. Corollary \ref{C} as a consequence of Theorem \ref{Thm3manif}.
\par\noindent
7. Terminology and basic facts about $3$-manifolds.
\par\noindent
8. Seifert foliations and pseudo-foliations.
\par\noindent
9. The annulus theorem and the JSJ decomposition.
\par\noindent
10. Digression on the terminology of the literature.
\end{abstract}

\maketitle

\section{\textbf{
Statement of the results
}}
\label{section1}

Consider a knot $K$ in $\mathbf S^3$.
Let $V_K$ be a tubular neighbourhood of $K$.
The \textbf{exterior} of $K$ is the closure $E_K$ of $\mathbf S^3 \smallsetminus V_K$,
and the \textbf{peripheral torus} is the common boundary 
$T_K = \partial V_K = \partial E_K$.
The \textbf{group of $K$} is the fundamental group $G_K = \pi_1(E_K)$,
and the \textbf{peripheral subgroup} is the image $P_K$ of $\pi_1(T_K)$ in $G_K$.
Recall that, by Dehn's Lemma, the map $\pi_1(T_K) \longrightarrow P_K$
is an isomorphism if and only if $K$ is non-trivial.

A subgroup $H$ of a group $G$ is \textbf{malnormal}
if $gHg^{-1} \cap G = \{e\}$ for all $g \in G$ with $g \notin H$;
basic facts on malnormal subgroups can be found
in our companion paper \cite{HaWeOs}.
The following question arose in discussions with Rinat Kashaev
(see also \cite{Kashaev} and  \cite{Kash--11});
we are grateful to him for this motivation.
\begin{itemize}
\item[ ]\emph{Given  $K$ as above,
when is  $P_K$ malnormal in $G_K$ ?}
\end{itemize}
The answer, our Corollary \ref{C}, happens to be a straightforward consequence 
of the following Theorem, from \cite{Simo--76};
the latter appears also as Lemma 1.1 in \cite{Whit--74} 
and Proposition 2 in \cite{Gram--91};
our proof, essentially self-contained, relies on Seifert foliations and pseudo-foliations.
Technical terms are defined below (see Sections
\ref{section2}, \ref{section3}, 
and the three appendices,
Sections \ref{section7},  \ref{SeifertF+psi}, and \ref{ann+JSJ}).

\begin{Thm}[\textbf{Reformulation of a result of Jonathan Simon}]
\label{laSimon}
Let $K$ be a non-trivial knot, $E_K$ its exterior, $T_K$ its boundary, 
and $\mu$ a meridian of $T_K$.
\par
Assume that there exists an annulus $A$ in $E_K$ attached to $T_K$
which is not boundary parallel.
Then the knot $K$ is 
\begin{itemize}
\item[(i)]
either a composite knot,
\item[(ii)]
or a torus knot,
\item[(iii)]
or a cable knot.
\end{itemize}
Moreover,
if $\eta$ denotes one of the components of $\partial A$:
in Case (i) $\eta$ is a meridian of $T_K$;
in Cases (ii) and (iii), the distance of $\mu$ and $\eta$
is  $\Delta(\mu,\eta) = 1$.
In particular, in all cases, $\Delta(\mu,\eta) \le 1$.
\par

Conversely, if $K$ is as in one of (i), (ii), and (iii),
then there exists an annulus $A$ in $E_K$ attached to $T_K$
which is not boundary parallel.
\end{Thm}

The ``converse part'' of the theorem is 
a rather straightforward consequence of the definitions,
see Section \ref{section3}.
As a corollary of Theorem \ref{laSimon} and of the annulus theorem:

\begin{Cor}
\label{C}
For a non-trivial knot $K$, the peripheral group $P_K$
is malnormal in $G_K$ if and only if $K$ is
neither a composite knot, nor a torus knot, nor a cable knot.
\end{Cor}

To view $P_K$ as a subgroup of $G_K$,
we need to choose a path 
from the base point in $E_K$ implicitely used to define $G_K$ 
to the base point in $T_K$ implicitely used to define $\pi_1(T_K)$,
so that $P_K$ is  a subgroup of $G_K$ defined up to conjugation only.
But the conclusion of Corollary \ref{C} 
makes sense since a subgroup and all its conjugates
are together malnormal or not.
A similar remark holds for the next theorem.

Theorem \ref{laSimon} suggests a result on malnormal peripheral subgroups
in a more general situation:

\begin{Thm}
\label{Thm3manif}
Let $M$ be a $3$-manifold which is 
compact, connected, orientable, and  irreducible. 
Assume moreover that the boundary $\partial M$ 
has at least one component, say $\partial_1 M$, which is a torus;
and that $M$ is
neither a solid torus $\mathbf S^1 \times \mathbf D^2$
nor a thickened torus $\mathbf T^2 \times [0,1]$.

Denote by $G$ the fundamental group of $M$,
by $P$ the peripheral subgroup of $G$ associated with $\partial_1 M$,
and by $V$ the connected component 
of the Jaco-Shalen-Johannson decomposition of $M$
which contains $\partial_1 M$.

Then $P$ is not malnormal in $G$ if and only if $V$ is a Seifert manifold.
\end{Thm}

Two observations are in order.
\par

In case $M$ is a solid torus or a thickened torus
(in both cases an irreducible Seifert manifold),
the peripheral group coincides with $\pi_1(M)$,
and thus is trivially malnormal in $\pi_1(M)$. 
\par

There is  a well-known fact
on 3-manifolds  which are compact, connected, orientable,
irreducible, and with non-empty boundary:
\emph{if one boundary component of such a manifold is a compressible torus,
then the manifold is a solid torus;}
for the convenience of the reader, we provide a proof of this
as Lemma  \ref{LemSolidTorus} below. 
Thus in the situation of Theorem \ref{Thm3manif},
$\partial_1 M$ is incompressible in $M$.

\medskip

By specialising to exteriors of links
(see the definitions recalled
at the end of Section \ref{SeifertF+psi}),
we could obtain the following corollary.
The notation we use for a link $L$, namely $V_L$, $E_L$ and $G_L$,
are defined in the same way as for knots.

\begin{Cor}
\label{Clink}
Let $L$ be a link in $\mathbf S^3$, with $r \ge 2$ components $L_1, \hdots, L_r$.
Assume that $L$ is unsplittable, and is not the Hopf link.
Denote by $G_L$ the group of $L$.
For $j \in \{1, \hdots, r\}$,
denote by $P_j$ the peripheral subgroup of $G_L$ which corresponds to $L_j$.
\par

Then $P_j$ is not malnormal in $G_L$ if and only if
\begin{itemize}
\item[$\circ$]
either $L_j$ is part of a (possibly satellised) torus sublink of $L$;
\item[$\circ$]
or $L_j$ is the outcome of a connected sum operation of links.
\end{itemize}
\end{Cor}

In this paper, 
we begin by giving a proof of Theorem \ref{laSimon}, 
following the method of \cite{Simo--76};
from this and the annulus theorem, Corollary \ref{C} follows.
Then, using more of the theory of $3$-manifolds,
we prove Theorem \ref{Thm3manif};
as a consequence, we obtain a second proof of Theorem \ref{laSimon}.
\par

More precisely, Section \ref{section2} contains general facts on annuli 
attached to boundary tori of $3$-manifolds.
Section \ref{section3} analyses exteriors of composite knots,
torus knots, and cable knots,
and thus establishes the converse (and easy) implication in Theorem \ref{laSimon}.
In Section \ref{section4}, we complete the proof of Theorem~\ref{laSimon};
the first few lines show also how Corollary \ref{C} follows from Theorem \ref{laSimon}.
Section \ref{section5} is a proof of Theorem \ref{Thm3manif},
and Section \ref{section6} shows how Corollary \ref{C} 
follows from Theorem \ref{Thm3manif}.
\par

We will not show how Corollary \ref{Clink} follows from Theorem \ref{Thm3manif},
for length reasons.
Indeed, if $T_1, \hdots, T_r$ denote 
the boundary components of the exterior $E_L$ of the link
$L = L_1 \sqcup \cdots \sqcup L_r$, 
a JSJ piece of $E_L$ can be adjacent to just one of the $T_j$
or to several of them, and many cases have to be treated separately,
so that there are (among other things) non-trivial combinatorial complications.
To avoid unreasonable length, 
we have chosen to leave the details  to the expert readers.
\par

As we have non-expert readers in mind,
we have written a rather long appendix, split in three parts.
In the first part, Section \ref{section7}, 
we recall various basic definitions on $3$-manifolds,
a theorem due to Alexander on complements of tori in $\mathbf S^3$,
and a re-embedding construction of Bonahon and Siebenmann
for submanifolds of $\mathbf S^3$.
Section \ref{SeifertF+psi} is about Seifert foliations 
and Seifert pseudo-foliations on $3$-manifolds.
The third part of the appendix, Section \ref{ann+JSJ},
is a reminder on the annulus theorem and the JSJ decomposition, 
needed for our proof of Theorem \ref{Thm3manif}.
The last section is a digression on the terminology.
\par

In a compagnion paper \cite{HaWeOs}, we collect basic facts
and (more or less) standard examples on malnormal subgroups
and on Frobenius groups of permutations.
\par

It is convenient to agree on the following \textbf{standing assumption}: 
\begin{center}
\emph{
all $3$-manifolds and surfaces below are assumed to be 
\\
compact, connected,  orientable,
and possibly with boundary,
}
\end{center}
unless either they are obviously not, 
such as links or boundaries, which need not be connected, 
or if it is explicitely stated otherwise,
as for the space of leaves of a Seifert foliation,
a surface which need not be orientable.
Moreover, maps, and in particular embeddings from one manifold into another,
are assumed to be \textbf{smooth}.

\medskip

We are grateful to Cameron Gordon and Roger Fenn
who provided a sketch of proof of Corollary \ref{C}
based on Theorem \ref{laSimon},
in discussions during the Symposium Michel Kervaire 
(Geneva, February 10-13, 2009).

\section{\textbf{
On annuli embedded in $3$-manifolds with torus boundaries
}}
\label{section2}

\subsection{Curves in tori and slopes.}
\label{slopes}

A simple closed curve in a surface is \textbf{essential}
if it is not homotopic to a point, 
equivalently if it does not bound an embedded disc.
Let $T$ be a $2$-dimensional torus;
a \textbf{slope} in $T$ is an isotopy class of essential simple closed curves.
These curves and slopes are non-oriented.
\par

The \textbf{distance}
$\Delta(s_1,s_2)$ of two slopes $s_1,s_2$ in $T$
is the absolute value of their intersection number,
namely
\begin{equation*}
\Delta(s_1,s_2) \, = \,
\min \left\{
\sharp (\sigma_1 \cap \sigma_2) \hskip.2cm \Big\vert \hskip.2cm
\aligned
&\text{$\sigma_j$ is a simple closed curve}
\\
&\text{representing $s_j$, $j=1,2$}
\endaligned
\right\} .
\end{equation*}
(This ``distance'' does not satisfy the triangle inequality,
but the terminology is however standard.)
Two slopes are isotopic if and only if their distance is zero;
two slopes (once oriented) define a basis of $H_1(T, \mathbf Z)$ 
if and only if their distance is one.
Observe that, if $\sigma_1, \sigma_2$ are two  essential simple closed curves in $T$
which are disjoint, and therefore isotopic, then
the closure of each connected component of
their complement $T \smallsetminus (\sigma_1 \cup \sigma_2)$
is an annulus embedded in $T$.
\par

For curves on tori and for slopes,
see \cite{Rolf--76} (in particular  Section 2.C)
and \cite{Boye--02}.
\par

In case a $2$-torus $T$ is given as the boundary of a solid torus,
a \textbf{meridian}  is
an essential simple closed curve $\mu$ on $T$
which bounds a disc in the solid torus,
and  a \textbf{parallel} is
an essential simple closed curve $\lambda$ on $T$
such that the homotopy classes of $\lambda$ and $\mu$, with orientations,
constitute a basis of $\pi_1(T) = H_1(T, \mathbf Z)$.

\subsection{Annuli attached to a torus component of the boundary.}
\label{annuliattached}

Let $M$ be a $3$-manifold with boundary,
such that one connected component of $\partial M$, say $T$, is a $2$-torus.
An \textbf{annulus in $M$ attached to $T$}
is an annulus $A$ which is properly embedded in $M$
and such that each component of $\partial A$ is an essential curve in $T$;
observe that these two components are disjoint, 
so that we have  a well-defined \textbf{slope of $A$ in $T$}.
\par

As a particular case of a definition from Subsection \ref{irred+para},
an annulus $A$ in $M$ attached to $T \subset \partial M$ 
is \textbf{boundary parallel}
if there exists a solid torus $U$ embedded in $M$ such that
\begin{itemize}
\item[(i)]
$A \subset \partial U$,
\item[(ii)]
$(\partial U \smallsetminus A) \subset T$,
\item[(iii)]
there exists a diffeomorphism\footnote{Note 
that $U$ and $A \times [0,1]$ are manifolds \emph{wich corners};
the notion of diffeomorphism has  to be adapated to this situation.}
$h : U \longrightarrow A \times [0,1]$
such that $h(A) = A \times \{0\}$;
\end{itemize}
in this case $A$ is said to be boundary parallel \textbf{through} $U$.
Observe that, in this situation, there is an
annulus $A_T$ 
embedded in $T$
such that $\partial U = A \cup_{\partial A} A_T$
(the notation $\cup_{\partial A}$ indicates that
$A \cap A_T = \partial A$).

Our next Subsections, \ref{subsection2.3}
to \ref{subsection3.3},
describe various examples in the particular situation of the exterior of a knot.

\subsection{Examples of boundary parallel annuli in knot exteriors.}
\label{subsection2.3}
Let $K$ be a knot and let $E_K, T_K = \partial E_K$ be as usual.
Any slope $s$ in $T_K$ can be the slope of a boundary parallel torus
attached to $T_K$ in $E_K$.
\par

Indeed, consider an annulus $A_T$ in $T_K$ bounded by 
two parallel essential simple closed curves in $T_K$ in the class $s$.
Push the interior of  $A_T$ slightly inside $E_K$ to obtain
an annulus $A$ in $M$ attached to $T_K$; 
there is a well-defined solid torus $U$, bounded by $A \cup A_T$,
such  that $A$ is boundary parallel through $U$.
\par

(This would essentially carry over to any boundary component of any $3$-manifold,
instead of just the torus $\partial E_K$.)
\par

On the contrary, Theorem \ref{laSimon} shows
that there are strong limitations on slopes  
which can be associated to non-boundary parallel annuli attached to $T_K$.
More generally there are strong limitations on slopes of incompressible surfaces,
see e.g. \cite{Boye--02}.

\subsection{Examples of annuli in the exterior of the trivial knot}
\label{subsection2.4}

Consider a non-trivial knot $J$, with $E_J, T_J$ as usual,
and a meridian $K \subset T_J$,
which is viewed as a trivial knot in $\mathbf S^3$,
with $E_K, T_K$ as usual.
Then $A := T_J \cap E_K$ is an annulus in $E_K$ attached to $T_K$.
It is not boundary parallel because $J$ is non-trivial.
The slope of $A$ in $T_K$ is a parallel.

\section{\textbf{
Examples of non boundary-parallel annuli 
\\
in knot exteriors
}}
\label{section3}

The three items below provide a proof of the converse part of Theorem \ref{laSimon}.

\subsection{Examples of  annuli in exteriors of composite knots}
\label{subsection3.1}

Consider a composite knot, namely a connected sum
$K = K_1 \sharp K_2$ of two non-trivial knots.
There is no loss of generality if we assume that $K$ is in $\mathbf R^3$,
intersecting $\mathbf R^2$ in exactly two points,
and that, if $H_1, H_2$ denote the two closed half-spaces 
bounded by $\mathbf R^2$, 
the knot $K_j$ is the union of $K \cap H_j$ 
with the straight segment in $\mathbf R^2$ joining the two points
of $K \cap \mathbf R^2$ (for $j = 1,2$).
Then 
\begin{equation*}
A \, = \,  
\big( \mathbf R^2 \smallsetminus (\mathbf R^2 \cap \overset{\circ}{V}_K) \big)
\cup \{\infty\}
\end{equation*}
is an annulus in $E_K$ attached to $T_K$.
\par
The slope of $A$ in $T_K$ is a meridian.
\par

For $j = 1,2$, denote by $W_j$ the closure of the complement in $H_j$
of $H_j \cap V_K$.
Observe that $W_j$ is diffeomorphic to the exterior of $K_j$
so that both $W_1$ and $W_2$ are knot exteriors.
\par

(Note: if $K_2$, say, were trivial, then $W_2$ would be a solid torus.)

\subsection{Examples of  annuli in exteriors of torus knots}
\label{subsection3.2}

Denote by $\mathbf S^1 \times \mathbf D^2$
the  solid torus standardly embedded in $\mathbf R^3$
and by $\mathbf T$ its boundary.
A \textbf{torus knot} is a knot isotopic to
an essential simple closed curve on the torus $\mathbf T$.
We agree here that
\begin{center}
\emph{the trivial knot is not a torus knot.}
\end{center}
\par

Let $K$ be a torus knot on $\mathbf T$, 
and let $V_K$ be a tubular neighbourhood of $K$
small enough for $V_K \cap \mathbf T$ to be a pair $(K_1, K_2)$
of curves on $\mathbf T$ which are disjoint and parallel to $K$
(the standardly embedded torus $\mathbf T$ 
should not be confused with
the torus $T_K = \partial V_K = \partial E_K$).
The complement $\mathbf T \smallsetminus (K_1 \cup K_2)$
has two connected components; 
let $A$ be the closure of the component
which does not contain $K$.
Then $A$ is an annulus in 
$E_K = (\mathbf R^3 \cup \{\infty\}) \smallsetminus \overset{\circ}{V}_K$
attached to the boundary $T_K$.
\par

The slope of $A$ in $T_K$ is  a parallel.
\par

The complement $E_K \smallsetminus A$ of $A$ has two connected components.
The bounded component is essentially the interior of the standard solid torus; 
more precisely it is the interior of this standard solid torus 
minus part of the ``small'' solid torus $V_K$;
thus, the closure of this bounded component is again a solid torus.
Similarly, the other component, together with the point at infinity of $\mathbf R^3$,
is a solid torus.

\subsection{Examples of  annuli in exteriors of cable knots}
\label{subsection3.3}

Consider on the one hand a non-trivial knot $K_c$
and a tubular neighbourhood $V_c$ of $K_c$
with its boundary $T_c = \partial V_c$.
Consider on the other hand
the standardly embedded solid torus 
$\mathbf S^1 \times \mathbf D^2$, 
a non-trivial torus knot $K_{pat}$ in $\partial (\mathbf S^1 \times \mathbf D^2)$,
and a homeomorphism $h : \mathbf S^1 \times \mathbf D^2 \longrightarrow V_c$.
Then, by definition,
the knot $K := h(K_{pat})$ is a \textbf{cable knot around $K_c$},
with \textbf{companion} $K_c$ and \textbf{pattern} $K_{pat}$.
We do assume that $K_c$ is non-trivial;
thus, in this paper,
\begin{center}
\emph{torus knots are not cable knots.}
\end{center}
Some authors (including \cite{Simo--76}) use the other convention,
and consequently state Theorem \ref{laSimon}
with two cases only.
%

\par

Let $A_{pat}$ be an annulus inside  $\partial (\mathbf S^1 \times \mathbf D^2)$
related to $K_{pat}$
as $A$ is related to $K$ in the previous Subsection \ref{subsection3.2}.
Then $A := h(A_{pat})$
is an annulus in $E_K$ attached to the boundary $T_K$. 
\par

The slope of $A$ in $T_K$ is again a parallel.
\par

The two components of $E_K \smallsetminus A$ are 
one a solid torus 
(which is a small perturbation of $V_c$),
and the other a knot exterior 
(a small perturbation of the exterior of $K_c$).


\section{\textbf{
Proof of Theorem \ref{laSimon}
}}
\label{section4}

We continue with the notation of Theorem \ref{laSimon}.
It is useful to consider a thickened torus 
\begin{equation*}
N_T \, := \,  T_K \times [0,\epsilon] \, \subset \, V_K
\hskip.3cm \text{with} \hskip.2cm N_T \cap E_K = T_K
= T_K \times \{0\},
\end{equation*}
as well as
a thickened annulus embedded in $E_K$
\begin{equation*}
N_A \, := \,  A \times [1,2] \, \subset \, E_K
\hskip.5cm \text{with} \hskip.2cm A = A \times \{\frac{3}{2}\} .
\end{equation*}
Define
\begin{itemize}
\item[-]
the shrinked neighbourbood $V_K^- := V_K \smallsetminus N _T$,
\item[-]
the enlarged exterior $E_K^+ :=  E_K \cup N_T = \overline{\mathbf S^3 \smallsetminus V_K^-}$,
\item[-]
and their common boundary
$T'_K := T_K \times \{\epsilon\} = \partial V_K^- = \partial E_K^+$,
which is $T_K$ slightly pushed inside $V_K$.
\end{itemize}
The union 
\begin{equation*}
N \,  := \,  N_T \cup N_A
\end{equation*}
is a manifold with boundary 
(indeed with corners).
Note that
\begin{equation*}
C_K \,  := \,  \partial A \times [1,2] = T_K \cap N_A = N_T \cap N_A = V_K \cap N_A
\end{equation*}
is the disjoint union of two annuli, each one being a neighbourhood in $T_K$
of a component of $\partial A$.
On $N$, there is a natural foliation by circles,
such that $\partial A \times \{1,2\}$ is the union of four particular leaves,
of which the isotopy class is a slope of $T_K$.
The manifolds $N$ and $N_A$ are irreducible 
since they are Seifert manifolds with boundary
(Proposition \ref{Seifertirr}).
\par

Let $\Theta$ denote the corresponding space of leaves.
Then $\Theta$ is homeomorphic to the union 
of an annulus (the space of leaves of $N_T$)
together with a thickened diameter (the space of leaves of $N_A$);
the four points which are common to the boundary of the thickened annulus
and the boundary of the thickened diameter
represent the four leaves in $\partial A \times \{1,2\}$.
Since $\Theta$ is orientable (Lemma \ref{seiforientable}), 
it is a planar surface with three boundary components (a ``pair of pants''),
and the manifold $N$ is diffeomorphic to a product:
\begin{equation*}
N \, \approx  \Theta \times \mathbf S^1 .
\end{equation*}
The boundary $\partial N$ is the union of three tori.
One is $T'_K$; 
we denote the two others by $T_1$ and $T_2$.
For $j \in \{1,2\}$, 
the torus $T_j$ separates $\mathbf S^3$ in two components;
we denote by $W_j$ the closure of the component contained in $E_K$.
Thus
\begin{equation*}
E_K^+ \, = \, N \cup W_1 \cup W_2 
\hskip.5cm \text{and} \hskip.5cm
E_K \, = \, N_A \cup W_1 \cup W_2
\end{equation*}
where the interiors on the right-hand sides are disjoint.
\par

By Alexander's Theorem \ref{ThmAlex}, 
each $W_j$ can be either a solid torus or a knot exterior,
so that there are three cases to consider:
\begin{itemize}
\item[(\ref{subsection4.1})]
both $W_1$ and $W_2$ are knot exteriors;
\item[(\ref{subsection4.2})]
both $W_1$ and $W_2$ are solid tori;
\item[(\ref{subsection4.3})]
$W_1$ is a solid torus and $W_2$ is a knot exterior.
\end{itemize}
We will see that
these three cases correspond respectively
to $K$ being a composite knot, a torus knot, and a cable knot.
Thus the proof below splits naturally in three cases;
it follows and extends the indications given by \cite{Simo--76}.

\subsection{Case in which both $W_1$ and $W_2$ are knot exteriors}
\label{subsection4.1}
We have to show that the slope of $A$ in $T_K$ is a meridian,
and it will follow that $K$ is a connected sum of two non-trivial knots.
Compare with Subsection \ref{subsection3.1}.

\par

Since $W_1$ is a knot exterior,
the manifold $\mathbf S^3 \smallsetminus \overset{\circ}{W}_1$
is  a solid torus, by Theorem \ref{ThmAlex}.
Thus $W_2$ is contained in the interior of a $3$-ball
which is  contained in $\mathbf S^3 \smallsetminus \overset{\circ}{W}_1$,
by Proposition \ref{BingMartin}; 
we denote by $\Sigma$ the boundary of this ball.
Since $\Sigma \cap (W_1 \cup W_2)$ is empty, we have
\begin{equation*}
\Sigma \, \subset \,  \text{interior of }  \left( V_K \cup N_A \right) 
\end{equation*}
and $\Sigma$ separates $W_1$ from $W_2$.
Since $V_K \cap N_A = C_K$, we have
\begin{equation*}
\Sigma \cap \overset{\circ}{C_K} \, \ne \, \emptyset .
\end{equation*}
Indeed, if this intersection were empty,
we would have either  $\Sigma \subset V_K$ or $\Sigma \subset N_A$,
and each of these inclusions would contradict the fact
that $\Sigma$ separates $W_1$ from $W_2$.
Moving slightly $\Sigma$ if necessary,
we can assume that 
the intersection $\Sigma \cap C_K$ is transverse,
so that this intersection, say $\Gamma$,  is a bunch of pairwise disjoint circles.
\par

We will  show how to modify $\Gamma$ (by modifying $\Sigma$),
in order to obtain a curve which is both a slope of $A$ in $T_K$
and a meridian of $T_K$.
\par

Consider a circle $\gamma$ of the bunch $\Gamma$
which is \textbf{innermost} in $\Sigma$, 
namely which bounds a disc, say $D_{\gamma} \subset \Sigma$,
such that $\overset{\circ}{D_{\gamma}} \cap \Gamma = \emptyset$.
There are two cases to consider, depending on this disc
being in $V_K$ or in $N_A$.
\par

Suppose first that $D_{\gamma} \subset N_A$.
Since the circle $\gamma$ is in one of the two annuli 
making up $C_K = \partial A \times [1,2]$, 
there are a priori two possibilities:
either it bounds a disc in this annulus, 
or it is parallel to the boundary of this annulus.
But the second case would mean that $\gamma$
defines the same slope of $T_K$ as $A$;
thus $\gamma$ would be essential in $T_K$;
since $\gamma$ bounds a disc in $E_K$,
the knot $K$ would be trivial,
in contradiction with our hypothesis.
Hence $\gamma$ bounds a disc
$D'_{\gamma} \subset \partial A \times \{j\}$, where $j=1$ or $j=2$.
The union $D_{\gamma} \cup D'_{\gamma}$ is a $2$-sphere
embedded in $N_A$;
\emph{since $N_A$ is irreducible,} 
this $2$-sphere  bounds a $3$-ball in $N_A$.
We can now isotope $D_{\gamma}$ through this $3$-ball
and then push it slightly outside $N_A$, to remove the intersection $\gamma$.
\par

Suppose now that $D_{\gamma} \subset V_K$.
As before, there are a priori two possibilites:
either $\gamma$ bounds a disc in $C_K$,
and we can modifiy the situation to one with one circle less,
or $\gamma$ is both a slope of $A$ and a meridian of $T_K$.
\par

Iterating the previous construction with an innermost circle
as often as necessary,
we obtain in all cases a curve which is both a slope of $A$
and a meridian of $T_K$.
This ends the proof  of Theorem \ref{laSimon} 
in case both $W_1$ and $W_2$ are knot exteriors.

\subsection{Case in which both $W_1$ and $W_2$ are solid tori}
\label{subsection4.2}
There are \emph{a priori} two subcases.
\par

Either the oriented foliation of $N$ introduced above 
extends to an oriented foliation by circles 
of $E_ K^+ = N \cup W_1 \cup W_2$. 
Then $K$ is a torus knot, by Proposition \ref{KnotSeifert}.
\par

Or the oriented foliation of $M$ 
does not extend to one of the $W_j$, say to $W_1$.
Then, by  Corollary \ref{uniquepseudo-leaf},  
this foliation extends to $W_2$, 
indeed to $\mathbf S^3 \smallsetminus \overset{\circ}{W_1}$,
with $K$ a regular leaf, and the knot $K$ is trivial.
Since the triviality of $K$ contradicts the hypotheses of Theorem \ref{laSimon},
this second subcase does not occur.

\subsection{Case in which  $W_1$ is a solid torus and $W_2$ a knot exterior}
\label{subsection4.3}
Let $\widehat N$ denote a submanifold of $\mathbf S^3$
obtained from $N$ by the Bonahon-Siebenmann re-embedding construction;
this amounts to replacing $W_2$ by a solid torus, that we denote by $U_2$
(see Proposition \ref{Bonahon-Siebenmann}).
Recall that $\widehat N$ is diffeomorphic to $N$,
and thus is given together with a Seifert foliation (indeed is a circle bundle).

\begin{Cla}
\label{ClaimSection4.3}
The foliation on $\widehat N$ extends to 
$\mathbf S^3 = \widehat N \cup V^-_K \cup W_1 \cup U_2$.
\end{Cla}

Let us admit the claim.
In the solid torus $W_1$,
the core must be an exceptional leaf,
otherwise the annulus $A$ would be 
boundary parallel through $W_1$.
Thus, if we consider the original embedding $N \subset \mathbf S^3$,
we see that $K$ is a cable around the knot whose exterior is $W_2$.

\medskip

\noindent
\emph{Proof of Claim \ref{ClaimSection4.3}.}
Since the complement of $\widehat N$ in $\mathbf S^3$ is a union of solid tori,
the Seifert fioliation on $\widehat N$ extends 
either as a Seifert foliation 
or as a pseudo-foliation, say $\mathcal F$, on $\mathbf S^3$.
By Corollary \ref{uniquepseudo-leaf},
it is enough to show that there cannot exist a pseudo-leaf
in any of $V^-_K$, $W_1$, $U_2$.
\par

Suppose that $V^-_K$ would contain a pseudo-leaf.
Then $\mathcal F$ would be the standard pseudo-foliation,
with unique pseudo-leaf inside $V^-_K$.
Hence the core of $W_1$ would be a regular leaf,
and $A$ would be  boundary parallel through\footnote{This
is abusive, since $W_1$ stands here for a slightly expanded solid torus $W^+_1$
made up of $W_1$ and the appropriate part of $N_A$.
But small lies can help the truth to be simpler.} 
the solid torus $W_1$. 
This would contradict the hypothesis on $A$, and is therefore impossible.
\par

The same argument shows that $U_2$ does not contain a pseudo-leaf.
\par

Suppose that $W_1$ contains a pseudo-leaf.
Then the knot $K$ is a regular leaf of the pseudo-foliation.
But regular leaves are all isotopic in $\mathbf S^3 \smallsetminus U_2$,
and they all bound discs.
Hence $K$ bounds a disc in $\mathbf S^3 \smallsetminus U_2$
which is a meridian disc in $V_K$.
This disc lives also in the original situation,
and this implies that the knot $K$ is trivial,
in contradiction with our hypothesis.
\hfill $\square$

\section{\textbf{
Proof of Corollary \ref{C} and of Theorem \ref{Thm3manif}
}}
\label{section5}

\medskip 

Consider a 3-manifold $M$ which satisfies
the hypothesis of Theorem~\ref{Thm3manif}, 
and assume moreover that $\partial_1 M \approx \mathbf T^2$ is incompressible
(see the comments which follow Theorem \ref{Thm3manif}).
\emph{We assume that $P \approx \pi_1(\partial_1 M)$ 
is not malnormal in $G = \pi_1(M)$,
and we have to show that the JSJ piece $V$ which contains $\partial_1 M$
is a Seifert manifold.}

\medskip

By assumption,
there exist  $p_0, p_1  \in P \smallsetminus \{1\}$
and $g \in G \smallsetminus P$
such that $g p_0 g^{-1} = p_1$.
The elements $p_0$ and $p_1$ are represented by loops in $\partial_1 M$
which are freely homotopic in $M$;
hence there exists a (possibly singular) map 
$\varphi : \mathbf A \longrightarrow M$
of which the image connects these two loops
($\mathbf A$ is the standard annulus $\mathbf S^1 \times [0,1]$).
\par

Let us check that $\varphi$ is essential
(compare with \cite{Simo--76}, Page 207).
On the one hand, one component of $\partial \mathbf A$
generates $\pi_1(\mathbf A) \approx \mathbf Z$;
its image by $\varphi$ is $p_0 \ne 1$ (or $p_1  \ne 1$) so that,
the group $P$ being torsion-free,
$\varphi$ induces an injection of $\pi_1(\mathbf A)$ into $P$,
and therefore also into $G$.
On the other hand, there is a spanning arc $\alpha$ in $\mathbf A$
which is mapped by $\varphi$ to $g$;
since $g \notin P$,
the restriction $\varphi \vert \alpha$ is not homotopic relative to its boundary
to an arc in $\partial M$.
\par


From the Annulus Theorem \ref{AnnulusThm},
there exists an  embedding 
$\psi :  \mathbf A \longrightarrow M$ 
with $\psi(\partial \mathbf A) \subset \partial_1 M$.
The annulus $\pi(\mathbf A)$
is not boundary parallel
(this is the meaning of $\psi$ being essential in Theorem \ref{AnnulusThm}).

\medskip
\noindent
\textbf{On the proof of Corolary \ref{C}.}
For a non-trivial knot $K$, 
the argument of the few lines above show that, if the peripheral subgroup $P_K$
is not malnormal in $G_K$, 
then there exists an annulus in $E_K$ attached to $T_K$
which is not boundary parallel.
Hence Corollary \ref{C} follows from Theorem \ref{laSimon}.
\par
\emph{Remark.}
The following is useless for our purpose but pleasant to know:
the peripheral subgroup $P_K$ of a knot group $G_K$
is maximal abelian in $G_K$ (a result of Noga, 
see Corollary~1 in \cite{Feus--70})
and of infinite index in $G_K$ (Theorem 10.6 in [Hemp--76]). 
\par

\medskip

We return to the situtation of Theorem \ref{Thm3manif}.
By Theorem \ref{JSJ}, it has a family $\mathcal T$ of tori
providing a JSJ decomposition in various pieces;
recall that $V$ denotes that piece which contains $\partial_1 M$.

\begin{Cla}
\label{Claim6}
There exists an  embedded essential annulus in $V$,
with at least one boundary component in $\partial_1 M$.
\end{Cla}

\noindent \emph{Proof.}
We know that there exists 
an  embedded essential annulus $\psi$ as above.
Without loss of generality, we can assume that
$\psi(\mathbf A)$ is transversal to $\mathcal T$,
so that the intersection $\psi(\mathbf A) \cap \mathcal T$ 
is a bunch $\mathcal B$ of circles.
Let us agree that such a circle $\beta$ is
\begin{itemize}
\item[$\bullet$]
of the \emph{first kind} if it bounds a disc in $\psi(\mathbf A)$,
\item[$\bullet$]
of the \emph{second kind} if it is boundary parallel in $\psi(\mathbf A)$.
\end{itemize}

First, we get rid of circles of the first kind,
by a classical argument.
As a preliminary observation, note that 
a circle from $\mathcal B$ contained inside a circle of the first kind
is also of the first kind.
Thus, if there are circles of the first kind, one may choose one of them, 
say $\beta$, which is innermost,
so that $\beta$ bounds a disc $\Delta$ in $\psi(\mathbf A)$ 
containing no element of $\mathcal B$ in its interior.
Denote  by $T_i$ the torus of the family $\mathcal T$
which contains $\beta$;
we have $\Delta \cap T_i = \partial \Delta$.
If $\beta$ was not contractible inside $T_i$,
the disc $\Delta$ would be a compressing disc for $T_i$,
and this is impossible since $T_i$ is incompressible.
Hence $\beta$ bounds a disc, say $\delta$,  in $T_i$.
The union $\delta \cup \Delta$ is a 2-sphere embedded in $M$;
\emph{since $M$ is irreducible,} it bounds a 3-ball.
We can first isotope $\Delta$ through this 3-ball 
and then push it slightly outside $T_i$
to remove the intersection $\beta$.
\par

Note that the irreducibility of $M$
has played a crucial role above.

Iterating this operation, we can obtain eventually
an annulus $\psi ' : \mathbf A \longrightarrow M$ 
embedded in $M$ with no circles of the first kind.
Note that $\psi ' (\mathbf A)$ and $\psi(\mathbf A)$ are isotopic in $M$,
so that $\psi '$ is also essential.
The intersection $\psi ' (\mathbf A) \cap \mathcal T$ is now 
a bunch $\mathcal B'$ of circles which are all of the second kind,
namely which are boundary parallel in $\psi '(\mathbf A)$.

The circles in $\mathcal B'$ decompose $\psi ' (\mathbf A)$ 
in a sequence of successive annuli.
If the first of them is essential, keep it and stop.
If it is inessential, then both its boundary components are in $\partial_1 M$,
and we can repeat the same argument with the second annulus.
After some time, we must encounter an essential annulus
with  one boundary component in $\partial_1 M$
and with empty intersection with $\mathcal T$,
therefore an essential annulus entirely contained in $V$.
\hfill $\square$

\begin{Cla}
\label{Claim7}
With the notation of the previous claim, the manifold $V$ is Seifert.
\end{Cla}

\noindent \emph{Proof.}
From the JSJ Decmposition Theorem,
we know that $V$ is Seifert or atoroidal.
But, if $V$ is atoroidal, it is also Seifert,
thanks to the following proposition
that we copy from Lemma 1.16 on Page 20 of \cite{Hatc--00}.
\hfill $\square$

\begin{Prop}
\label{Prop8}
Let $X$ be a compact,  irreducible, and atoroidal $3$-manifold.
Assume that $X$ contains an incompressible, $\partial$-incompressible annulus
meeting only torus components of $\partial M$.
\par
Then $X$ is a Seifert manifold.
\end{Prop}

This ends the proof of the implication 
``$P$ not malnormal $\Longrightarrow$ $V$ Seifert''
of Theorem~\ref{Thm3manif}.

\medskip
\noindent
\textbf{Proof of the converse statement in Theorem \ref{Thm3manif}.}
We assume that $V$ is a Seifert manifold.
It is an irreducible manifold:
one reason is noted as Comment (i) after Theorem \ref{JSJ},
another one is that $V$ has a boundary 
(Proposition \ref{Seifertirr}).
\par

We have the inclusions $\partial_1 M \subset V \subset M$
and the corresponding group homomorphisms
$\mathbf Z^2 \approx \pi_1(\partial_1 M)  
\longrightarrow \pi_1(V)  \longrightarrow \pi_1(M) = G$.
As already noted in Section~\ref{section1},
we can assume that $\partial_1 M$ 
is incompressible in $M$, 
so that the inclusion $\pi_1(\partial_1 M) \longrightarrow G$
is an isomorphism onto the peripheral subgroup $P$ of $G$;
a fortiori, the inclusion $\pi_1(\partial_1 M) \longrightarrow \pi_1(V)$
is an injection.
The homomorphism $\pi_1(V)  \longrightarrow \pi_1(M)$
is also an injection (Remark (iv) after Theorem \ref{JSJ}).
Let us moreover remark that  $\pi_1(\partial_1 M)$
is a proper subgroup of $\pi_1(V)$;
otherwise, because of standard facts on fundamental groups
of Seifert manifolds with boundaries (see Chapter 12 in \cite{Hemp--76}),
$V$ would be a thickened torus, and this has been ruled out. 
Summing up: $\pi_1(\partial_1 M) \approx \mathbf Z^2$
is a proper subgroup of $\pi_1(V)$.
\par

Since $V$ is a Seifert manifold, the torsion-free group $\pi_1(V)$ 
has an infinite cyclic normal subgroup,
which is generated by the homotopy class of a regular fibre.
By Proposition 2.viii of \cite{HaWeOs}, it follows that $\pi_1(V)$
does not have any non-trivial malnormal subgroup.
A fortiori, $P$ is not malnormal in~$G$.

\section{\textbf{
Corollary \ref{C} as a consequence of Theorem \ref{Thm3manif}
}}
\label{section6}

Corollary \ref{C}  is a consequence of Theorem \ref{Thm3manif} 
and of various facts on Seifert foliations (see Section \ref{SeifertF+psi}) 
which are summed up in the following proposition.

\begin{Prop}
\label{SeifertJSJext}
Let $K$ be a non-trivial knot,
$E_K$ its exterior, 
$T_K$ its boundary,
and $V$ the component of the JSJ decomposition of $E_K$ containing $T_K$.
\par
If $V$ is a Seifert manifold, then $K$ is
either a composite knot, or a torus knot, or a cable knot.
\end{Prop}

\noindent
\emph{Note.}
One proof would be to show that, if $V$ is Seifert,
then there exists an essential embedded annulus in $V$.
We will proceed differently, without using Theorem \ref{laSimon}.

There are statements similar to our proposition
in \cite{JaSh--79L}, Lemma VI.3.4,
and \cite{Joha--79}, Lemma 14.8.
The proof below is somewhat different, 
as it uses Seifert foliations and pseudo-foliations.

\medskip

\noindent
\emph{Proof of Proposition \ref{SeifertJSJext}.}
Denote by $T_1, \hdots, T_r$ the connected components of $\partial V$
distinct from $T_K$ (possibly $r=0$, in case of $K$ a torus knot, 
since then $E_K$ is Seifert foliated).
For $j \in \{1, \hdots, r\}$, 
let $W_j$ denote the closure of the connected component of 
$\mathbf S^3 \smallsetminus T_j$ which does not contain $V$.
If $W_j$ was a solid torus, $T_j$ would be compressible in $E_K$;
but this would contradict the incompressibility 
of the JSJ tori in the decomposition of $E_K$.
Hence, by Theorem \ref{ThmAlex}, $W_j$ is  a knot exterior.
\par

We consider the Bonahon-Siebenmann re-embedding 
$\widehat V$ of $V$ in $\mathbf S^3$, 
summed up below in Proposition \ref{Bonahon-Siebenmann}.
This construction amounts to replace each $W_j$ by a solid torus $U_j$.
Then $\widehat V$ can be seen as the exterior $E_L$
of a link $L := \widehat K \sqcup L_1 \sqcup \cdots \sqcup L_r$
with $r+1$ components.
\par

Proposition \ref{SeifertJSJext} is now a consequence of the following claim.
\hfill $\square$

\begin{Cla}
\label{ClaimSeifertJSJext}
(i) If the Seifert foliation on $V$ does not extend to $V_K$
(or, equivalently, if the Seifert foliation on $\widehat V$
does not extend  to $V_K$),
then $r \ge 2$ and $K$ is the connected sum of $r$ prime knots.
\par

(ii) If the Seifert foliation on $V$ extends to $V_K$, then $r \le 1$.
If $r=0$, then $K$ is a torus knot, and if $r=1$ then $K$ is a cable knot.
\end{Cla}

\medskip

\noindent
\emph{Proof.}
(i) By hypothesis, there exists a foliation on $\widehat V$ 
which does not extend to $V_K$.
By Proposition \ref{LinkSeifert}, 
this foliation does extend as a pseudo-foliation on $\mathbf S^3$,
with $\widehat K$ as a pseudo-leaf.
By the uniqueness result for pseudo-foliations of $\mathbf S^3$, 
Corollary \ref{uniquepseudo-leaf},
the cores of the solid tori $U_i$ are regular leaves.
It follows from the description of composite knots \`a la Schubert
(Subsection \ref{laSchubert}) that $K$ is a connected sum of $r$ prime knots. 
\par

(ii) By hypothesis, there exists a foliation on $V$ which extends to $M = V_K \cup V$;
observe that the manifold $M$ is irreducible 
(being Seifert and with boundary, Proposition \ref{Seifertirr}).
To show that $r \le 1$, 
we proceed by contradiction and assume for some time that $r \ge 2$. 
\par

Consider the torus component $T_1$ of $\partial M$.
Since $M$ is irreducible and is not a solid torus,
$T_1$ is incompressible in $M$ (Lemma \ref{LemSolidTorus}).
Since $T_1$ bounds on the other side the cube with a knotted hole $W_1$,
it is incompressible in $W_1$.
Thus, van Kampen's theorem shows that $\pi_1(T_1)$ is a subgroup
of $\pi_1(M \cup W_1)$.
By a similar argument, $\pi_1(T_1)$ is still a subgroup of $\pi_1(M \cup W_1 \cup \cdots \cup W_r)$.
But this is impossible since
$M \cup W_1 \cup \cdots \cup W_r = \mathbf S^3$.
\par

If $r=0$, then $K$ is a torus knot, since it is a leaf of a Seifert foliation
of $\mathbf S^3$.
\par

If $r=1$, the same incompressibility argument as above
shows that $M$ is a solid torus
(otherwise $T_1$ would be incompressible in $\mathbf S^3$).
We know the classification of Seifert foliations on a solid torus:
the space of leaves is a disc, and the number $s$ of exceptional leaves is at most~$1$.
If one had $s=0$, the manifold $V$ would be a thickened torus,
and this is impossible because $T_1$ cannot be boundary parallel
in a JSJ decomposition.
Hence $s=1$, and the knot $K$ is not an exceptional leaf
(this would again imply that $V$ is a thickened torus).
Hence $K$ is a regular leaf of the Seifert foliation on the solid torus $M$,
and this foliation has one exceptional leaf (the core of $M$).
Thus $K$ is a torus knot in $M$, not isotopic inside $M$ to the core of $M$.
It follows that this torus knot is satellised around the knot $K_1$
of which the exterior is $W_1$. This is exactly the cable situation.
\hfill $\square$

\bigskip

\begin{center}
\textbf{
Appendix: a reminder of some three-dimensional topology
}
\end{center}

\section{\textbf{
Terminology and basic facts about $3$-manifolds
}}
\label{section7}

This section is a reminder on
some terminology for $3$-manifolds,
and classical results that we have used in 
Sections \ref{section4}, \ref{section5} and \ref{section6}
(Alexander, Dehn, Seifert, Waldhausen, Bing-Martin, Bonahon-Siebenmann).
Recall the \textbf{standing assumption}
agreed upon in Section \ref{section1}:
\begin{center}
\emph{
all $3$-manifolds and surfaces below
are assumed
to be compact, connected, orientable,
and possibly with boundary,
}
\end{center}
unless a few exceptions which are either obvious or explicitely stated as such.
Also, all maps are assumed to be smooth.
\par

A map $\varphi$ from a manifold $N$ to a manifold $M$
is \textbf{proper} if $\varphi^{-1}(\partial M) = \partial N$. 
A manifold $S$ is \textbf{properly embedded} in a manifold $M$
if it is embedded and if 
$\partial S = S \cap \partial M$.

\subsection{Irreducibility and parallelism}
\label{irred+para}

A $3$-mani\-fold $M$ is \textbf{irreducible} if 
any embedded $2$-sphere bounds a $3$-ball.
For example, $\mathbf S^3$ is irreducible;
indeed, it is a \textbf{theorem of Alexander}
that any $2$-sphere embedded in $\mathbf S^3$
bounds two $3$-balls
(see e.g.~Theorem 1.1 in \cite{Hatc--00}, Page 1).
%
The only irreducible $3$-manifold 
which has a $2$-sphere in its boundary is the $3$-ball.
For the importance of irreducibility hypothesis above,
see for example near the end of Subsection \ref{subsection4.1},
or the proof of Claim \ref{Claim6}.
\par

Let $M$ be a manifold of dimension $m$.
Let $S_0, S_1$ be two manifolds of the same dimension $n < m$,
with $S_0$ properly embedded in $M$
and $S_1$ either properly embedded in $M$ or embedded in $\partial M$.
Then $S_0$ and $S_1$ are \textbf{parallel} 
if there exists an embedding of a thickened manifold
\begin{equation*}
\psi \, : \,  S \times [0,1] \longrightarrow M
\end{equation*}
such that 
\begin{itemize}
\item[(i)]
$\psi(S \times \{0\}) = S_0$ and $\psi(S \times \{1\}) = S_1$, 
\item[(ii)]
$\psi(\partial S \times [0,1]) \subset \partial M$.
\end{itemize}
If  $\partial M \ne \emptyset$,
a manifold $S_0$ properly embedded in $M$
is \textbf{boundary parallel}, or \textbf{$\partial$-parallel},
if there exists a manifold $S_1$ embedded in $\partial M$
such that $S_0$ and $S_1$ are parallel.
\par

Consider for example the case with $M = A$ an annulus and $S$ of dimension $1$.
There are three isotopy classes of properly embedded arcs in an annulus $A$:
one class with the two ends of the arc in one component of $\partial A$,
one class with the two ends of the arc in the other component of $\partial A$,
these are boundary parallel, 
and the class of the so-called \textbf{spanning arcs} 
with one end in each component of $\partial A$,
equivalently with $A \smallsetminus \alpha$ connected.
\par

Recall that a simple closed curve in a surface is \textbf{essential}
if it is not homotopic to a point, 
equivalently if it does not bound an embedded disc.
A circle embedded in $A$  is essential 
if and only if it is boundary parallel, and is then a \textbf{core} of $A$.
Note that a core of $A$ and a spanning arc of $A$, appropriately oriented,
have intersection number $+1$.

\subsection{Incompressible and $\partial$-incompressible surfaces}
\label{incom+deltaincom}

Let $S$ be a surface properly embedded in a $3$-manifold $M$
and $\gamma$ a simple closed curve in the interior of $S$.
A \textbf{compressing disc} for $\gamma$ is a disc $D$ embedded in $M$
such that $\partial D = \gamma$ and $\partial D = D \cap S$.
The surface $S$ is \textbf{incompressible} if,
for any simple closed curve $\gamma$ in the interior of $S$
which has a compressing disc $D$,
there exist a disc $D'$ in $S$ such that $\partial D' = \gamma$
(equivalently: $\gamma$ is null-homotopic in $S$).
Note that our definition is different from that of \cite{Hemp--76} for
properly embedded surfaces which are discs or spheres;
for us, these are \emph{always} incompressible.
\par

Mutatis mutandis, this definition of ``incompressible'' applies
to boundary components of $M$.
\par
A non-connected surface (for example $\partial M$ in some situations)
is incompressible if each of its connected components is so
(see e.g.~Section 1.2 in \cite{Hatc--00}).
\par

A connected surface $S$ properly embedded in $M$ 
is incompressible if and only if
the induced homomorphism of groups
$\pi_1(S) \longrightarrow \pi_1(M)$ is injective.
This follows from  \textbf{Dehn's Lemma and the loop theorem};
see for example Corollary 3.3 in \cite{Hatc--00}.
(It is important here that $S$ is two-sided,
but this follows from our \emph{standing assumption,}
according to which both $S$ and $M$ are orientable.)
\par

For example, the boundary $T_K$ 
of a non-trivial knot exterior $E_K$ is incompressible, 
and the boundary of a handlebody of genus $g \ge 1$ is compressible.
\par

Given a (not necessarily connected) surface $S$ properly
embedded in a 3-manifold $M$,
the manifold $M^*_S$ obtained from $M$ by \textbf{splitting $M$ along $S$}
is the complement in $M$ of a regular open neighbourhood of $S$
(observe that $S$ is two-sided, being orientable in an orientable manifold).
We quote now Theorem 1.8 in \cite{Wa--67ab}.

\begin{Prop}[\textbf{Waldhausen}]
\label{splittingirreducibleWald}
Let  $M$ and $S$ be as above; assume that $S$ is incompressible.
Then the connected components of $M^*_S$ are irreducible 
if and only if $M$ is irreducible.
\end{Prop}

Let $S$ be a surface with boundary $\partial S \ne \emptyset$
properly embedded in a $3$-manifold $M$ with boundary $\partial M \ne \emptyset$.
For an arc $\alpha$ properly embedded in $S$,
a \textbf{compressing disc} is a disc $D$ embedded in $M$ with:
\begin{itemize}
\item[$\circ$]
$\alpha = D \cap S$, 
\item[$\circ$]
$\beta := D \cap \partial M$ is an arc in $\partial D$,
\item[$\circ$]
$\partial D = \alpha \cup \beta$ 
and $\partial \alpha = \partial \beta =  \{\text{two points in } \partial D \}$.
\end{itemize}
(Observe that such a $D$ is never properly embedded in $M$
since the interior of $\alpha$ is disjoint from $\partial M$.)
The surface $S$ is \textbf{$\partial$-incompressible} if,
for any arc $\alpha$ properly embedded in $S$ with $D$ as above,
$\alpha$ is boundary parallel in $S$.

\begin{Prop}
\label{onannuli}
Let $M$ be an irreducible $3$-manifold.
Assume that the boundary of $M$ has some torus compoments;
let $\partial_T M$ denote the union of these.
\par

Let $S$ be a surface properly embedded and incompressible in $M$,
with $\emptyset \ne \partial S \subset \partial_T M$.
Then either $S$ is $\partial$-incompressible
or $S$ is a boundary parallel annulus.
\par

In particular, if $\partial M$ is a union of tori,
an annulus properly embedded and incompressible in $M$
is either $\partial$-incompressible or boundary parallel.
\end{Prop}

\noindent
\emph{Proof.} We refer to Lemma 1.10 of \cite{Hatc--00}.
\hfill $\square$

\medskip

\noindent
\emph{Remark.} If $\partial_T M$ is compressible, 
it follows from Lemma \ref{LemSolidTorus} below that $M$ is a solid torus.
It is known that, in a solid torus, an incompressible surface which is not a disc
is necessarily 
an annulus parallel to the boundary
(Lemma 2.3 in \cite{Wa--67ab}).

\subsection{Complements of tori in the $3$-sphere}
\label{comptorisphere}

Let $T$ be a torus embedded in $\mathbf S^3$.
By  Poincar\'e-Alexander duality, 
the complement $\mathbf S^3 \smallsetminus T$
has two connected components,
and their closures $U_1, U_2$
have $T$ as a common boundary.
By the theorem of Alexander recalled in Subsection \ref{irred+para},
the manifolds $U_1$ and $U_2$ are irreducible.

The following Theorem \ref{ThmAlex} 
is also due to Alexander (see \cite{Hatc--00}, Page 11).
The proof below (unlike that of Alexander!)
uses Dehn's Lemma.
Our preliminary Lemma \ref{LemSolidTorus} 
is well-known to specialists.

\begin{Lem}
\label{LemSolidTorus} 
Let $M$ be an irreducible  3-manifold;
assume that the boundary $\partial M$ 
has a component $\partial_1 M$ which is a compressible torus.
\par

Then $M$ is a solid torus; in particular, $\partial M$ is connected.
\end{Lem}

\noindent \emph{Proof.}
Let $D$ be a compressing disc for $\partial_1 M$
and let $E$ be a small open tubular neighbourhood of $D$.
Let $M^*_D = M \smallsetminus E$ be the result of splitting $M$ along $D$.
By construction, the boundary $\partial M^*_D$ contains a 2-sphere,
consisting of  ``most of'' $\partial_1 M \cup \partial E$.
By the irreducibility assumption,  this 2-sphere (viewed now in $M$)
bounds a 3-ball $B$ in $M$.
Then $V \Doteq B \cup E$ is a solid torus, 
because it is obtained by attaching $E$ 
along $\partial B$ as a $1$-handle, 
and because it is orientable.
\par

This solid torus is closed in $M$ since it is compact.
It is also open by Brouwer's theorem of invariance of domain
(see the remark below). It follows that $V = M$.
\hfill $\square$


\medskip

\emph{Remark on Brouwer's Theorem.}
The following result is (a restatement of what is)
found in books, see e.g. Proposition 7.4 in \cite{Dold--72}):
\emph{an injective continuous mapping 
$g : N_1 \longrightarrow N_2$
between two manifolds $N_1, N_2$, 
of the same dimension and without boundary, is open.} 
Let now $M_1, M_2$ be manifolds of the same dimension, 
with boundary, and let $\partial ' M_2$ be the union of some
of the connected components of $\partial M_2$. Then:
\emph{an injective continuous mapping 
$f : M_1 \longrightarrow M_2$
such that $f(\partial M_1) = \partial ' M_2$  is open.}
This is a straightforward consequence of the previous statement,
applied to the natural map $g$ induced by $f$,
with domain the double $N_1 = M_1 \cup_{\partial M_1} M_1$ of $M_1$
and target the interior of the double $N_2 = M_2 \cup_{\partial ' M_2} M_2$.

\medskip  

The \textbf{core} of a solid torus $U$ embedded in a $3$-manifold $M$
is $h(\mathbf S^1 \times \{0\})$, 
where $U$ is the image of an embedding
$h : \mathbf S^1 \times \mathbf D^2 \longrightarrow  M$
of the standard solid torus
(the core is well-defined up to isotopy).

\begin{Thm}[\textbf{Alexander}]
\label{ThmAlex}
Let $T$ be a torus embedded in $\mathbf S^3$
and let  $U_1, U_2$ be the closures of the connected components
of $\mathbf S^3 \smallsetminus T$.
\par
At least one of them, say $U_1$, is a solid torus, say with core $C_1$,
so that $U_2$ is the exterior of the core   $C_1$.
The curve $C_1$ is unknotted in $\mathbf S^3$
if and only if $U_2$ is also a solid torus.
\end{Thm}

\noindent \emph{Proof.}
If $T$ were incompressible in both $U_1$ and $U_2$,
the group $\pi_1(T)$ would inject in $\pi_1(U_1)$ and $\pi_1(U_2)$,
by Dehn's lemma.
Since $\mathbf S^3 = U_1 \cup_T U_2$,
it would also inject in the amalgamted sum
\begin{equation*}
\pi_1(\mathbf S^3) 
\, = \, 
\pi_1(U_1) \ast_{\pi_1(T)} \pi_1(U_2),
\end{equation*}
by the Seifert--van Kampen theorem,
and this is absurd.
Upon exchanging $U_1$ and $U_2$,
we can therefore assume that $T$ is compressible in $U_1$.
Lemma \ref{LemSolidTorus} implies that $U_1$ is
a solid torus.
\par

If $U_2$ is also a solid torus,
then $T$ is unknotted, by definition. \hfill $\square$

\medskip

Alexander's theorem is strongly used in the following construction,
that we propose to call
the \textbf{Bonahon-Siebenmann's re-embedding construction}.
In \cite{BoSi} (see the first 10 lines or so in their Section 2.2), 
this is called a \emph{splitting.}
On page 326 of \cite{Budn--06}
and with the notation of our Proposition \ref{Bonahon-Siebenmann}, 
the embedding of $\widehat Z$ in $\mathbf S^3$ 
is called the \emph{untwisted re-embedding.}
\par

Let $Z$ be a $3$-manifold embedded in $\mathbf S^3$,
with boundary  a non-empty disjoint union of tori 
$\partial Z = T_1 \sqcup \hdots \sqcup T_r$.
For $j \in \{1, \hdots, r\}$, denote by $W_j$ the closure of
the connected component of $\mathbf S^3 \smallsetminus T_j$
which does not contain $Z$. 
Assume that the notation is such that $W_1, \hdots, W_\ell$ are knot exteriors
and $W_{\ell+1}, \hdots, W_r$ solid tori, for some $\ell$ with $0 \le \ell \le r$.
The purpose of the construction is to obtain a situation with $\ell = 0$,
namely with $Z$ re-embedded as the exterior of an appropriate link in the $3$-sphere.
\par

For $j \in \{1, \hdots, \ell\}$, denote by $\mu_j$ a meridian and by $\pi_j$ a parallel
of $T_j$ viewed as the boundary
of the solid torus $\overline{\mathbf S^3 \smallsetminus W_j}$;
orient these so that they become a basis of $H_1(T_j, \mathbf Z)$.
Let $U_j$ denote a solid torus, with boundary endowed with 
an oriented meridian $\mu'_j$ and an oriented parallel $\pi'_j$.
Define inductively a sequence of manifolds
$M_0 = \mathbf S^3, M_1, \hdots, M_\ell$;
for $j \in \{1, \hdots, \ell\}$,
the manifold $M_j$ is obtained by 
gluing $U_j$ to the closure of $M_{j-1} \smallsetminus W_j$,
in such a way that $\mu_j$ is glued to $\pi'_j$ and $\pi_j$ to $\mu'_j$.
Observe that the construction provides an embedding of  $Z$ in $M_j$,
and that the components of the complement of the image of $Z$ in $M_j$
can be naturally identified with
$U_1, \hdots, U_j, W_{j+1}, \hdots, W_r$.
Since $M_j$ has a Heegaard decomposition of genus one\footnote{A
manifold $N$ which has a Heegaard decomposition of genus one
is diffeomorphic to either the $3$-sphere, 
or $\mathbf S^1 \times \mathbf S^2$, 
or a lens space. 
Indeed, the median torus of such a decomposition 
is the boundary of two solid tori, and has therefore two meridians $\mu, \mu'$.
If $\delta = \Delta(\mu, \mu')$, with the notation of Subsection \ref{slopes},
then $\pi_1(N) \approx \mathbf Z /\delta \mathbf Z$,
and $N$ is $\mathbf S^1 \times \mathbf S^2$, or $\mathbf S^3$, or a lens space,
if the value of $\delta$ is $0$, or $1$, or $\ge 2$.
See Chapter 2 of \cite{Hemp--76}.}
and has the homology of the $3$-sphere,
$M_j$ is diffeomorphic to $\mathbf S^3$.
In particular, $M_\ell$ is diffeomorphic to $\mathbf S^3$,
and we denote it by $\mathbf S^3$ again;
we denote by $\widehat Z$ 
the image of the embedding of $Z$ in this ``new'' $\mathbf S^3$.
\par
For reference, we state
the result of this construction as: 

\begin{Prop}
\label{Bonahon-Siebenmann}
Let $Z$ be a 
$3$-manifold embedded in $\mathbf S^3$,
with boundary  a non-empty disjoint union of tori.
\par

There exists a submanifold $\widehat Z$ of $\mathbf S^3$
which is diffeomorphic to $Z$
and which is the exterior of a link in $\mathbf S^3$.
\end{Prop}

Using  automorphisms of the solid tori we could show that,
given two results $\widehat Z, \widehat Z' \subset \mathbf S^3$ of the construction,
there exists a diffeomorphism $h$ of $\mathbf S^3$
such that $h(\widehat Z) = \widehat Z'$.

\medskip

The following result is due to Bing and Martin.
A manifold with boundary homeomorphic 
to the exterior of a non-trivial knot in $\mathbf S^3$
is picturesquely called in \cite{BiMa--71}
a \textbf{cube with a knotted hole}; 
we also use \textbf{knot exterior}.
We insist that a knot exterior is  the exterior of a \emph{non-trivial} knot;
observe that Bing-Martin's result 
does not carry over to the exterior of a trivial knot.

\begin{Prop}[\textbf{Bing-Martin}]
\label{BingMartin}
Consider a solid torus $U$ in $\mathbf S^3$
and a knot exterior $W$ contained in the interior of $U$.
\par

Then there exists a $3$-ball $B$ in the interior of $U$
such that $W \subset \overset{\circ}{B}$. 
\end{Prop}

\noindent
\emph{Proof.}
Since $W$ is the exterior of a non-trivial knot,
the inclusion of $\partial W$ in $W$ induces an injection
of $\pi_1(\partial W) \approx \mathbf Z^2$ into $\pi_1(W)$,
by Dehn's lemma (the proof of Bing and Martin \emph{does not} use Dehn's lemma). 
Since 
\begin{equation*}
\mathbf Z \approx \pi_1(U) \, \approx \,  
\pi_1(U \smallsetminus \overset{\circ}{W}) \ast_{\pi_1(\partial W)} \pi_1(W) ,
\end{equation*}
by the Seifert--van Kampen theorem,
this implies that $\pi_1(\partial W)$ does not inject in
$\pi_1(U \smallsetminus \overset{\circ}{W})$.
Hence, by Dehn's lemma, there exists a compressing disc $D$
for $\partial W$ in $U \smallsetminus \overset{\circ}{W}$.
The union of $W$ and of a thickening of this disc $D$
is the $3$-ball we are looking for.
\hfill $\square$

\section{\textbf{
Seifert foliations and pseudo-foliations
}}
\label{SeifertF+psi}

\subsection{Seifert foliations}
\label{Seifertfoliations}

In this paper, a foliation always means a 
\textbf{foliation by circles of a $3$-manifold} $M$,
namely a partition of $M$ in circles with the usual regularity hypothesis.
A foliation is \textbf{oriented} if all its leaves are coherently oriented.
\par

For example, fixed-point free actions of the rotation group $SO(2)$
on $3$-manifold provide oriented foliations (see Proposition \ref{bbEpstein}).
\par

Standard actions of $SO(2)$ on the solid torus provide important examples.
More precisely (we follow Page 299 of \cite{OrRa--68}),
parametrise the solid torus $\mathbf S^1 \times \mathbf D^2$
by $(e^{i\psi}, \rho e^{i\theta})$, with
$0 \le \rho \le 1$ and $0 \le \psi, \theta < 2\pi$.
Given two coprime integers $\mu,\nu$ with $0 \le \nu \le \mu$, 
the corresponding \textbf{standard linear action} of
$SO(2) = \{z \in \mathbf C \hskip.1cm \vert \hskip.1cm \vert z \vert = 1 \}$
on the solid torus is defined by
\begin{equation*}
SO(2) \times \mathbf S^1 \times \mathbf D^2 
\, \longrightarrow \,
\mathbf S^1 \times \mathbf D^2 ,
\hskip.5cm
(z, e^{i\psi}, \rho e^{i\theta}) 
\, \longmapsto \,
( z^{\mu} e^{i\psi}, z^{\nu} \rho e^{i\theta} ).
\end{equation*}
This action is always effective\footnote{Recall that
an action of a group $G$ on a set $X$ is \emph{effective}
(or \emph{faithful}) if,
for any $g \in G$, $g \ne e$, there exists $x \in X$ with $gx \ne x$.
}.
It is free if and only if $\mu = 1$ (this implies $\nu = 0$ or $\nu = 1$);
in this case, the orbits are fibers of a product fibration
$\mathbf S^1 \times \mathbf D^2  \longrightarrow \mathbf D^2$.
If $\mu > 1$ (this implies $1 \le \nu \le \mu - 1$),
the action is free on the complement of the core of equation $\rho = 0$;
this core is the \emph{exceptional orbit},
and the other orbits are the \emph{regular} ones.
\par

For $\mu = 0$ and $\nu = 1$,
the same formulas define an action
with regular orbits in meridian discs,
and with the core as fixed point set.
See Subsection \ref{Seifertpseudo} below.
\par

A \textbf{Seifert foliation} is a circle foliation 
such that each leaf has a neighbourhood which is a union of leaves
and which is isomorphic to a standard foliated solid torus.
Each neighbourhood isomorphic to a standard foliated solid torus
with $\mu \ge 2$ contains an \textbf{exceptional leaf},
and all other leaves are  \textbf{regular leaves}. 
Many authors use  ``Seifert fibration''  for our ``Seifert foliations'',
but our terminology is motivated by the fact that these
in general \emph{are not} circle bundles,
and because ``fibration'' is already used in many other situations.
\par

Let $M$ be a $3$-manifold given with a Seifert foliation.
Let $\mathcal B$ be the corresponding \textbf{space of leaves}, 
viewed here as a surface with boundary, possibly non-orientable 
(we \emph{do not} view $\mathcal B$ as an orbifold),
and let $p : M \longrightarrow \mathcal B$ denote the quotient map.
Since $M$ is orientable, $\mathcal B$ is orientable 
if and only if the Seifert foliation is orientable.
We mark the points in $\mathcal B$ which correspond to exceptional leaves in $M$
and we denote by $\check{\mathcal B}$ the surface obtained from $\mathcal B$
by removing disjoint small discs around the marked points.
It is standard that $\mathcal B$, with appropriate decorations,
provides a complete description of $M$ and its Seifert foliation
(see e.g. Theorem 2.4 in \cite{Bona--02}).
\par

The following proposition is not deep, but useful.
It  shows an equivalence between Seifert foliations and leaves on the one hand,
and fixed point free $SO(2)$-actions and orbits on the other hand.

\begin{Prop}
\label{bbEpstein}
The leaves of an oriented Seifert foliation on a $3$-manifold $M$
are the orbits of a fixed point free action of $SO(2)$ on $M$, 
and conversely.
\end{Prop}

\noindent \emph{Proof.}
Any fixed point free $SO(2)$-action gives rise to an oriented Seifert foliation;
this is an immediate consequence of the existence of a slice for the action
(see \cite{Kosz--53}, 
or the middle of Page 304 in \cite{OrRa--68}).
\par

For the converse implication, we can quote \cite{Raym--68} or \cite{OrRa--68},
who prove that any Seifert data can be realised by a $SO(2)$-action.
Alternatively, see \cite{Epst--72}, Page 80, for a more direct approach.
\hfill $\square$

\medskip

The next theorem, from \cite{Epst--72}, is  much deeper
than Proposition \ref{bbEpstein}.
It is remarkable enough to be stated here, 
even if we do not use it elsewhere in this paper.
(Recall that our $3$-manifolds are compact, connected, and orientable.)

\begin{Thm}[\textbf{Epstein}]
\label{Epstein}
The leaves of an oriented circle foliation on a $3$-manifold $M$
are the orbits of a fixed point free action of $SO(2)$ on $M$.
\end{Thm}

The following result is standard.
For the first claim, see  e.g. Lemma VI.7 in \cite{Jaco--80}
or Proposition 1.12 in \cite{Hatc--00}.
For the second claim, see Corollary 3.2 \cite{Scot--83}.

\begin{Prop}
\label{Seifertirr}
A Seifert manifold is either irreducible, 
or  $\mathbf S^1 \times \mathbf S^2$,
or the connected sum of two projective spaces.
\par

In particular, a Seifert manifold with boundary is irreducible.
\par

Moreover, if $M$ is a Seifert manifold with boundary,
either $\partial M$ is incompressible or $M$ is a solid torus.
\end{Prop}

Finally, let us quote
a proposition which shows that
there are plenty of incompressible annuli and tori in Seifert manifolds.
For the proof, see Page 127 in \cite{JaSh--79L}.

\begin{Prop}
\label{examplesverticala+t}
Let $M$ be a Seifert manifold, and let $p : M \longrightarrow \mathcal B$
denote the projection on its space of leaves.
Let $\alpha$ be a properly embedded arc in $\mathcal B$
which avoids the marked points; 
set $A = p^{-1}(\alpha)$.
Let $\gamma$ be a simple closed curve in $\mathcal B$
which avoids the marked points
and which is orientation-preserving; 
set $T = p^{-1}(\gamma)$.
Then :
\par

(i) $A$ is a properly embedded annulus in $M$ which is incompressible.
\par
(ii) $A$ is boundary parallel in $M$ if and only if 
$\alpha$ is boundary parallel in $\mathcal B$.
\par
(iii) $T$ is a properly embedded torus in $M$; 
 it is compressible in $M$ if and only if 
$\gamma$ bounds a disc in $\mathcal B$ 
which contains at most one marked point.
\end{Prop}

These are examples of so-called \emph{vertical}
annuli and tori in Seifert manifolds.
For a precise description of annuli and tori in $3$-manifolds,
see Theorems 3.9 and 3.5 in \cite{Bona--02}, respectively.

\subsection{Seifert pseudo-foliations, general facts}
\label{Seifertpseudo}

Consider the standard solid torus $\mathbf U := \mathbf S^1 \times \mathbf D^2$
and its core $C = \{0\} \times \mathbf S^1$.
The \textbf{standard Seifert pseudo-foliation} of $\mathbf U$ 
is the partition $\mathcal F_0$ of $\mathbf U$ in the points of the circle $C$
and the circles  $\gamma_\rho \times \{z\}$,
where $\gamma_\rho$ is a circle of  radius $\rho > 0$ in $\mathbf D^2$
centred around the origin and $z$ is a point of $\mathbf S^1$;
we call $C$ the \textbf{pseudo-leaf} 
and the other circles the \textbf{leaves}
of  $\mathcal F_0$.
\par

More generally, let us define a \textbf{Seifert pseudo-foliation}
of a $3$-manifold $M$ 
to be a partition of $M$ in circles and points,
which restricts to a Seifert foliation 
outside a finite disjoint union of solid tori $V_1, \hdots, V_r$, 
and to standard Seifert pseudo-foliations on these solid tori.
Pseudoleaves and leaves of such a pseudo-foliation
are defined naturally; 
leaves can be either regular or exceptional, as in Seifert foliations.
We insist that we assume $r \ge 1$; in other words:
\begin{center}
\emph{
it is part of the definition of a Seifert pseudo-foliation 
that it contains at least one pseudo-leaf.}
\end{center}
By definition, a Seifert pseudo-foliation of $M$ as above 
restricts to a Seifert foliation on
$M \smallsetminus (\overset{\circ}{V_1} \sqcup \cdots \sqcup \overset{\circ}{V_r})$.
From now on, we will write \textbf{pseudo-foliation} 
instead of Seifert pseudo-foliation.

\par

For example, 
and as a consequence of the classification
%
%
of circle foliations on the $2$-torus,
any circle foliation $\mathcal F_0$ on the boundary of a solid torus 
extends to a Seifert foliation or a pseudo-foliation $\mathcal F$ on the solid torus itself.
More precisely,
$\mathcal F$ is a pseudo-foliation 
if the leaves of $\mathcal F_0$ are meridians,
and $\mathcal F$ is a Seifert foliation in all the other cases.
Hence, 
any Seifert foliation or pseudo-foliation on a submanifold of $\mathbf S^3$
with boundary a union of tori
extends to a foliation or a pseudo-foliation to $\mathbf S^3$ itself.

A Seifert  pseudo-foliation is \textbf{oriented} 
if its restriction to the complement of the pseudo-leaves is oriented.
\par

Proposition \ref{bbEpstein} has an analogue for pseudo-foliations:

\begin{Prop}
\label{pseudobbEpstein}
The leaves of an oriented pseudo-foliation on a $3$-manifold $M$
are the orbits of an effective action of $SO(2)$ on $M$ with fixed points, 
and conversely.
\end{Prop}

\noindent
\emph{Proof}.
That any effective action of $SO(2)$ with fixed points
gives rise to a pseudo-foliation is again a straightforward consequence
of the slice theorem (see the proof of Proposition \ref{bbEpstein}).
For the converse, we believe
that Epstein's proof carries over.
Alternatively, see Theorem 1 in [OrRa--68] and Corollary 2b, with $t = 0$,
of [Raym--68].
\hfill $\square$

\medskip

The \textbf{space of orbits} $\mathcal B$ of a pseudo-foliation of a manifold $M$
is again a surface (possibly non orientable) with boundary;
we use also ``space of leaves'' instead of ``space of orbits'',
even though this is abusive, 
since the restriction of the projection $M \longrightarrow \mathcal B$
to each pseudo-leaf  is a bijection.
There are again marked points in $\mathcal B$,
corresponding to exceptional  leaves in $M$
and we define $\check{\mathcal B}$ as above;
also, it is again true (as in \ref{Seifertfoliations})
that $\mathcal B$ together with appropriate decorations 
provides a complete description of $M$ and its pseudo-foliation.
Here are two basic examples.
\par

(i)
For the standard pseudo-foliation of the solid torus $\mathbf U$,
the space of orbits is an annulus.
One component of its boundary is 
the space of orbits of the $2$-torus $\partial \mathbf U$,
the other component corresponds to the pseudo-leaf,
which is the core of~$\mathbf U$,
and there are no marked points.
The canonical projection $p : \mathbf U \longrightarrow \mathcal B$
restricts to a bijection from the core of $\mathbf U$
onto one component of the boundary of the annulus $\mathcal B$.
\par
Conversely, let $M$ be a manifold with a pseudo-foliation 
such that the corresponding space $\mathcal B$ is an annulus
with one boundary component being the space of orbits of $\partial M$,
and without marked points;
then $M$ is a solid torus with the standard pseudo-foliation.
\par

(ii) 
On $\mathbf S^3 = \mathbf R^3 \cup \{\infty\}$, 
the standard action of $SO(2)$ by rotations around an axis
defines a pseudo-foliation with one pseudo-leaf.
The space of orbits  $\mathcal B$ is the $2$-disc $\mathbf D^2$,
and the boundary of the disc represents the fixed points of the action.
Conversely, if a pseudo-foliation on $M$ gives rise to such a $2$-disc,
then $M$ is the $3$-sphere, 
the leaves are the circular orbits of the standard action of $SO(2)$,
and the pseudo-leaf is the circle of fixed points.

\medskip 

To simplify the discussion, 
we assume from now on that $M$ has no boundary.
In particular, boundary components of $\mathcal B$
are in bijection with pseudo-leaves of $M$.

\subsection{Seifert pseudo-foliations on closed $3$-manifolds}
\label{Seifertpseudoclosed}
To emphasise the difference between foliations and pseudo-foliations,
we state as a lemma an observation due to Waldhausen
(see
Page 90 of \cite{Wa--67ab}, the discussion of Condition 6.2.4).

\begin{Lem}[\textbf{Waldhausen}]
\label{Waldhausen}
Let $\mathcal F$ be a pseudo-foliation on a closed $3$-mani\-fold $M$,
let $p : M \longrightarrow \mathcal B$ 
be the projection on the space of orbits,
and let $\beta$ be an arc properly embedded in $\mathcal B$
which avoids the marked points.
\par 

Then $p^{-1}(\beta)$ is a $2$-sphere embedded in $M$.
\end{Lem}

\noindent
\emph{Proof.}
Remove a little interval at each extremity of $\beta$.
Let $\beta^*$ denote the closure of the complement of these intervals.
Then $p^{-1}(\beta^*)$ is an annulus.
The inverse image of a little interval is a disc,
indeed a meridian disc in the tubular neighbourhood of the pseudo-leaf.
Thus $p^{-1}(\beta)$ is an annulus with a disc 
glued on each of its boundary components; hence $p^{-1}(\beta)$ is a $2$-sphere.
\hfill $\square$

\medskip

Easy and standard topological considerations show that
Lemma \ref{Waldhausen} has the following consequences.

\begin{Prop}
\label{ConsequencesWaldhausen}
Let the notation be as in the previous lemma.
\begin{itemize}
\item[(i)]
The sphere $p^{-1}(\beta)$ separates $M$
if and only if $\beta$ separates $\mathcal B$,
thus producing a connected sum decomposition of $M$.
\item[(ii)]
If $\beta$ does not separate $\mathcal B$, then the sphere $p^{-1}(\beta)$
produces a factor $\mathbf S^1 \times \mathbf S^2$ in $M$.
\item[(iii)]
If $\mathcal B$ is a disc with one marked point,
then $M$ is a lens space.
\item[(iv)]
The sphere $p^{-1}(\beta)$ bounds a $3$-ball
if and only if $\beta$ is boundary parallel in $\check{\mathcal B}$.
\end{itemize}
\end{Prop}

\noindent
\emph{Comments.}
Claim (i) is straightforward.
Claim (ii) holds by classical arguments (Lemma 3.8, Page 27, in \cite{Hemp--76}).
For (iii), see Page 301 of \cite{OrRa--68}.
Claim (iv) follows from the previous ones.
\hfill $\square$

\medskip

The following result can either be easily deduced from the considerations above,
or recovered as a special case of a result of Orlik and Raymond,
written in terms of $SO(2)$-actions with fixed points
(Page 299 of \cite{OrRa--68}).

\begin{Prop}
\label{Mwithpseudo}
Let $M$ be a closed orientable manifold which admits an orientable pseudo-foliation.
\par
Then $M$ is either a $3$-sphere, or $\mathbf S^1 \times \mathbf S^2$,
or a lens space, or a connected sum of these.
\end{Prop}

\begin{Cor}
\label{uniquepseudo-leaf}
Let $M$ be a homology $3$-sphere; assume that $M$
is furnished with a pseudo-foliation.
\par
Then $M = \mathbf S^3$ is the standard sphere,
and the pseudo-foliation, up to diffeomorphism,
is given by the standard action of the rotation group $SO(2)$
around an axis of the sphere $\mathbf S^3 = \mathbf R^3 \cup \{\infty\}$,
as in Example (ii) of the end of \ref{Seifertpseudo}.
\par
In particular, the pseudo-foliation has exactly one pseudo-leaf.
\end{Cor}

\noindent 
\emph{Proof.}
Let $\mathcal F$ be a pseudo-foliation on $M$;
by our next Lemma \ref{seiforientable}, $\mathcal F$ is orientable. 
\par
If the space of leaves $\mathcal B$ of $\mathcal F$
either had at least two boundary components
or had genus $\ge 1$,
there would exist in $\mathcal B$ an arc $\beta$ as in Lemma \ref{Waldhausen},
and Proposition \ref{ConsequencesWaldhausen}.ii 
would imply that $\mathbf Z = H_1(\mathbf S^1 \times \mathbf S^2, \mathbf Z)$ 
is a direct factor of $H_1(M, \mathbf Z)$;
this is impossible because $M$ is a homology sphere.
Hence $\mathcal B$ is a disc.
\par

If this disc had just one marked point,
$M$ would be a lens space by Proposition \ref{ConsequencesWaldhausen}.iii,
and this is again impossible since $M$ is a homology sphere.
If this disc had two or more marked points,
$M$ would be a connected sum of lens spaces,
and this is equally impossible.
Hence $\mathcal B$ is a disc without marked points.
\par

This implies that $M = \mathbf S^3$,
and that $SO(2)$ acts as stated in the corollary.
\hfill $\square$

\medskip

A Seifert foliation (or pseudo-foliation!)
on the $3$-sphere is necessarily orientable.
More generally:

\begin{Lem}
\label{seiforientable}
Let $Z$ be a $3$-manifold embedded in a homology $3$-sphere~$M$.
\par
If $Z$ has a Seifert foliation with space of leaves~$\mathcal B$, 
then $\mathcal B$ is orientable.
\end{Lem}

\noindent
\emph{Proof.}
Let us first recall that, in a homology $3$-sphere $M$,
any embedded surface $S$ without boundary is orientable.
Indeed, using homology and cohomology with coefficients 
the field $\mathbf{F}_2$ with two elements 
(so that $H^2(S) \approx \mathbf{F}_2$
for any closed surface $S$, orientable or not),
Alexander duality
\begin{equation*}
H^2(S) \, \approx \, H_0(M \smallsetminus S) / H_0(\text{point})
\, \approx \, \mathbf F_2
\end{equation*}
shows that $S$ is two-sided, and therefore orientable.
\par

To prove the lemma, it is enough to show that
any simple closed curve $\gamma$ in $\mathcal B$
which avoids the marked points is two-sided.
Given such a $\gamma$, 
consider the surface $S := p^{-1}(\gamma)$.
Since $S$ is orientable in an orientable manifold,
$S$ is two-sided; it follows that $\gamma$ is two-sided.
\par
[Note that, since $S$ is orientable and is a circle bundle over a circle,
$S$ is a torus].
\hfill $\square$

\medskip

We end this subsection by translating
Propositions  \ref{Seifertirr} and  \ref{Mwithpseudo}
in terms of $SO(2)$-actions:

\begin{Cor}
Let $M$ be a closed orientable $3$-manifold.
\begin{itemize}
\item[(i)]
If $M$ affords a fixed-point free $SO(2)$-action,
then $M$ is either irreducible, or $\mathbf S^1 \times \mathbf S^2$,
or the connected sum of two projective spaces.
\item[(ii)]
If $M$ affords a non-trivial $SO(2)$-action with fixed points,
then $M$ is  $\mathbf S^3$, or $\mathbf S^1 \times \mathbf S^2$,
or a lens space, or a connected sum of these.
\item[(iii)]
If $M$ affords $SO(2)$-actions of the two kinds, with fixed points and without,
then $M$ is  $\mathbf S^3$, or $\mathbf S^1 \times \mathbf S^2$, 
or a lens space.
\end{itemize}
\end{Cor}

Thus, the list of $SO(2)$-actions with fixed points if far more restricted
than the list of actions without fixed points.

\medskip

It is a natural temptation to hope for a general theory
which would encompass Seifert foliations and pseudo-foliations;
but we should not give in this, as Waldhausen has warned us
(see the \emph{Bemerkung} on Page 91 of \cite{Wa--67ab}).
Indeed, Seifert manifolds are irreducible
(up to a small number of exceptions, see Proposition \ref{Seifertirr});
irreducibility is a crucial ingredient of their theory and classification.
On the contrary, ``most'' pseudo-foliated manifolds 
are reducible (Proposition \ref{Mwithpseudo});
as a consequence, this ``general theory'' would be worthless.

\subsection{Seifert manifolds embedded in $\mathbf S^3$}
We begin by stating a standard characterisation of torus knots,
used above in Subsection \ref{subsection4.2}.

\begin{Prop}[\textbf{Seifert foliations on knot exteriors}]
\label{KnotSeifert}
A  knot $K$ such that $E_K$ carries a Seifert foliation 
is a torus knot or the trivial knot.
\end{Prop}

\noindent
\emph{Proof.} 
There are two cases to distinguish.
\par
(i) Suppose that the foliation extends to $V_K$, 
providing a Seifert foliation of $\mathbf S^3$.
By Seifert's classification of the Seifert foliations on the $3$-sphere \cite{Seif--33},
$K$ is either a torus knot or the trivial knot.
\par
(ii) Suppose that the foliation, say $\mathcal F$,  does not extend to $V_K$.
Then the induced foliation on $T_K$ is necessarily a foliation by meridians, 
so that $\mathcal F$ extends to a pseudo-foliation $\mathcal F'$ on $\mathbf S^3$.
By Corollary \ref{uniquepseudo-leaf}, $\mathcal F'$ has a unique pseudo-leaf,
which is $K$ and which is a trivial knot.
\hfill $\square$ 

\medskip


Proposition \ref{KnotSeifert} suggests to distinguish two types of links, as follows
(it is a result from \cite{BuMu--70}).

\begin{Prop}[\textbf{Seifert foliations on link exteriors}]
\label{LinkSeifert}
Let $L = L_1 \sqcup \cdots \sqcup L_r$ be a link in $\mathbf S^3$,
with  $r \ge 2$ components.
Assume that $E_L$ admits a Seifert foliation $\mathcal F$.
\begin{itemize}
\item[(i)]
If $\mathcal F$ extends to a foliation of the $3$-sphere,
then $L$ is obtained by selecting $r$ leaves of a Seifert foliation of $\mathbf S^3$.
Links of this type are called \textbf{torus links}.
\item[(ii)]
If $\mathcal F$ does not extend as a foliation,
then it extends to a pseudo-foliation of the $3$-sphere, 
necessarily with a unique pseudo-leaf
which is a component of $L$, say $L_1$.
The other components $L_2, \hdots, L_r$ of $L$
are meridians around the pseudo-leaf $L_1$. 
\end{itemize}
\end{Prop}

\noindent
\emph{Proof.} This is a straightforward consequence 
of Corollary \ref{uniquepseudo-leaf}. 
\hfill $\square$

\medskip

Let us finally define  terms which enter our Corollary \ref{Clink}.
A link $L$ is \textbf{unsplittable} if,
for any 2-sphere  $S$ in $\mathbf S^3$ disjoint from $L$, 
all components of the link $L$
are in the same connected component of 
$\mathbf S^3 \smallsetminus S$.
A \textbf{torus sublink} 
is a part $L_{i_1} \sqcup \cdots \sqcup L_{i_s}$
of a torus link $L_1 \sqcup \cdots \sqcup L_r$,
for some subsequence of indices with
$1 \le i_1 < \cdots < i_s \le r$.
A \textbf{connected sum operation} on two links
$L = L_1 \sqcup \cdots \sqcup L_r$ and $L' = L'_1 \sqcup \cdots \sqcup L'_s$
consists in connecting \emph{just one} component $L_j$ of $L$
with one component $L'_k$ of $L'$.
A \textbf{satellised sublink} of $L$ is a link obtained from $L$
by replacing some of the $L_j$'s
by satellites
of these $L_j$'s.

\subsection{Composite knots \`a la Schubert}
\label{laSchubert}
We revisit Schubert's description of composite knots
in terms of pseudo-foliations and the re-embedding construction.
Consider a solid torus $\mathbf U$ embedded in an unknotted way in $\mathbf S^3$,
an integer $r \ge 2$, 
and $r$ disjoint  meridian discs $D_1, \hdots, D_r$ in $\mathbf U$;
thicken these discs a little bit; the thickened discs
separate $\mathbf U$ in $r$ closed $3$-balls $B_1, \hdots, B_r$.
\par

Let $K$ be a knot embedded in the interior of $\mathbf U$.
We assume that $K$ intersects each disc $D_i$ transversely
in exactly one point.
Hence $K$ runs across each thickened disc
in a little unknotted arc,
and $A_i := K \cap B_i$ is a properly embedded arc in $B_i$, 
for $i = 1, \hdots, r$.
We assume moreover that the arc $A_i$ is knotted in $B_i$.
Denote by $K_i$ the knot obtained from $A_i$ by adding
an unknotted arc outside $B_i$.
Then, by construction-definition, the knot $K$
is the connected sum $K_1 \sharp \cdots \sharp K_r$
of the knots $K_1, \hdots, K_r$.
\par

Denote by $V_K$ a thin tubular neighbourhood of $K$ in $\mathbf U$,
and set $T_K = \partial V_K$.
Consider a little collar of $T_K$ inside $V_K$, 
denote by $T'_K$ the component of its boundary which is inside $V_K$, 
and by $V^-_K$ the smaller tubular neighbourhood of $K$
with boundary $T'_K$.

For $i \in \{1, \hdots, r\}$, 
let $W_i$ be the closure of $B_i \smallsetminus (V_K \cap B_i)$.
Thus $W_i$ is a cube with a knotted hole;
this hole is indeed knotted, since $A_i$ is knotted by hypothesis.
Note that $W_i$ is the exterior of the knot $K_i$ defined above.
Set $T_i := \partial W_i$, which is a $2$-torus.
\par

On the one hand, define
\begin{equation*}
\Sigma \, := \, \mathbf S^3 \, \smallsetminus \,
\left( V^-_K \sqcup W_1 \sqcup \hdots \sqcup W_r \right)
\end{equation*}
(the $\sqcup$ indicate disjoint unions),
and observe that $\partial \Sigma = T'_K \sqcup T_1 \sqcup \cdots \sqcup T_r$.
On the other hand, consider the link
$L = L_0 \sqcup L_1 \sqcup \cdots \sqcup L_r$ in $\mathbf S^3$ 
obtained from the standard pseudo-foliation
by selecting the pseudo-leaf $L_0$ and $m$ regular leaves $L_1, \hdots, L_r$.
The Bonahon-Siebenmann re-embedding 
$\widehat \Sigma$ of $\Sigma$ in $\mathbf S^3$
is obtained by replacing each cube with a knotted hole $W_i$
by a cube with an \emph{unknotted hole},
namely by a solid torus, say $U_i$.
Thus $\widehat \Sigma$ is the exterior of the link $L$ defined just above.
\par

From that description, we see that the link exterior $\widehat \Sigma$
is diffeomorphic to a product $F \times \mathbf S^1$,
where $F$ is a planar surface with $r+1$ boundary components
$\partial_0 F, \hdots, \partial_r F$.
The product $\partial_i F \times \mathbf S^1$ 
is the boundary of a little tubular neighbourhood of $L_i$.
Hence $\Sigma$ is also Seifert foliated and diffeomorphic to $F \times \mathbf S^1$.
Since the foliation on $\widehat \Sigma$ 
does not extend to a neigbourhood of $L_0 = \widehat K$,
the foliation $\Sigma$ does not extend to a neibhbourhood of $K$. 
In both foliations (of $\Sigma$ and of $\widehat \Sigma$),
the leaves on the boundary of a tubular neighbourhood are meridians.

To obtain the composite knot $K$ from the link $L$,
we have just to replace each cube with an unknotted hole $U_i$
by $W_i$. After this replacement, 
the link component $L_0 = \widehat K$ is changed into $K$.

\emph{Remark on the necessity of the condition $r \ge 2$:}
In the connected sum point of view,
we wish the sum to be non-trivial, namely involving at least one non-trivial factor.
In the JSJ point of view, if the Seifert foliation on $V$ does not extend to $V_K$
and if $r=1$, then the JSJ torus $T_1$ is boundary parallel (parallel to $T_K$),
and this contradicts the JSJ conditions.

\section{\textbf{
The annulus theorem and the JSJ decomposition
}}
\label{ann+JSJ}

\subsection{The annulus theorem}
\label{annulustheorem}

A first major ingredient of our proof of Theorem \ref{Thm3manif}
is the annulus theorem,
announced together with the torus theorem by Waldhausen \cite{Wald--69};
a detailed proof of the annulus theorem
was given in \cite{CaFe--76}.
Before stating the theorem, we recall some terminology.
\par

We denote by $\mathbf A$ the standard annulus $\mathbf S^1 \times [0,1]$,
and by $\partial_0 \mathbf A = \mathbf S^1 \times \{0\}$,
$\partial_1 \mathbf A = \mathbf S^1 \times \{1\}$
the two components of its boundary $\partial \mathbf A$.
Recall that \emph{spanning arcs} in $\mathbf A$ 
have been defined in \ref{irred+para}. 
A proper map from the standard annulus to a $3$-manifold $M$, 
say $\varphi : (\mathbf A, \partial \mathbf A) \longrightarrow (M, \partial M)$,
is \textbf{essential} if the induced homomorphism
$\pi_1(\mathbf A) \longrightarrow \pi_1(M)$ is injective 
and if, for a spanning arc $\alpha$ in $\mathbf A$, 
the restriction $\varphi \vert \alpha$
is not homotopic rel its boundary to an arc in $\partial M$.

 \begin{Thm}[\textbf{Annulus Theorem}]
\label{AnnulusThm}
Let $M$ be a compact orientable $3$-manifold 
and let $\varphi : \mathbf A \longrightarrow M$
be an essential map from the annulus  to~$M$.

Then there exists an essential embedding $\psi : \mathbf A \longrightarrow M$
such that, for $i=0$ and $i=1$, the image $\psi(\partial_i \mathbf A)$ 
lies in the same connected component of $\partial M$ as $\varphi(\partial_i \mathbf A)$.
\end{Thm}

%

\noindent
\emph{On the proof.} The last part of our formulation, from ``such that'',
is~not explicit in Theorem 3 of \cite{CaFe--76},
but it follows from their proof. \hfill $\square$

\subsection{The JSJ decomposition}
\label{JSJdecomposition}

The second major ingredient in our proof of Theorem \ref{Thm3manif}
is the \textbf{JSJ  decomposition of $3$-manifolds},
as stated below.
Traditionally, JSJ refers explicitely to Jaco-Shalen-Johan\-n\-son,
and also implicitely to Waldhausen \cite{Wald--69}.


A 3-manifold $M$ is \textbf{atoroidal} (other authors use ``\textbf{simple}'')
if any incompressible torus in $M$ is boundary parallel 
(as defined in \ref{irred+para}).
Recall that $M^*_{\mathcal T}$ denotes the manifold obtained
by splitting $M$ along $\mathcal T$.

\begin{Thm}[\textbf{JSJ Decomposition}]
\label{JSJ} 
Let $M$ be an irreducible manifold with empty or incompressible boundary.
\par

In the interior of $M$, there exists  
a family $\mathcal T = \{T_1, \hdots, T_r\}$ of disjoint tori
which are incompressible  and not parallel to components of $\partial M$, 
with the following properties:
\begin{itemize}
\item[(i)]
each component of $M^*_{\mathcal T}$ is either a Seifert manifold
or atoroidal;
\item[(ii)]
the family $\mathcal T$ is minimal among those
which have Property (i).
\end{itemize}
Moreover, 
such a family $\mathcal T$ is unique up to ambient isotopy.
\end{Thm}

\noindent \emph{Proof.} 
See (among others) 
Theorem 1.9 in \cite{Hatc--00}
and Page 169 of \cite{JaSh--79L}.
See also the comments in the much shorter \cite{JaSh--79C},
or Theorem 3.4 of \cite{Bona--02}.
\hfill $\square$
\medskip

The connected components of $M^*_{\mathcal T}$ are the
\textbf{pieces of the JSJ-decomposition},
and $\mathcal T$ is the \textbf{characteristic torus family}.
Observe that, by Condition (ii), 
no piece can be a thickened torus,
unless $M$ itself is a thickened torus
(in which case $\mathcal T$ is the empty family,
and $M$ has a unique component, itself).

\medskip

\noindent
\emph{Comments on the statement of Theorems \ref{JSJ}.}
(i) By Proposition \ref{splittingirreducibleWald},
each piece is irreducible.
\par
(ii)
A $3$-manifold can be both atoroidal and Seifert;
examples include the solid torus, the thickened torus,  
exteriors of torus knots,
and some manifolds without boundary.
The complete list, which is rather short, can be found in \cite{JaSh--79L}
(Page 129, and Lemma IV.2.G).
Thus, the ``or'' in (i) above is not exclusive.
\par

(iii)
Let $V$ be a piece.
Since the $T_j$'s are incompressible,
the inclusion $V \subset M$
induces an injection of $\pi_1(V)$ into $\pi_1(M)$.
This follows from an appropriate version of the
Seifert--van Kampen theorem.

\medskip

\noindent
\emph{Beyond Theorems \ref{AnnulusThm}, \ref{JSJ}, and \ref{Thm3manif}.}
The Cannon-Feustel theorem \ref{AnnulusThm}
can be thought of as the original annulus theorem.
There exists a stronger result, due to Johannson 
(the ``enclosing theorem'' of \cite{Joha--79})
and Jaco-Shalen (the ``mapping theorem'' of \cite{JaSh--79L},
see their remark in the middle of Page 56);
see also Page 173 of \cite{Jaco--80}.
This requires more general JSJ-like decompositions,
with a characteristic surface $\mathcal A$ (rather than $\mathcal T$)
composed of \emph{annuli and tori},
and with pieces which can be 
\emph{ atoroidal, or Seifert manifolds or $\mathbf I$-bundles}
($\mathbf I$ is the unit interval).
The existence and unicity of such an $\mathcal A$
is closely related to the strong result just alluded to,
and to a general homotopy annulus theorem
for essential mapping from a Seifert manifold to a $3$-manifold.
For a statement which is both precise and readable,
we refer to Theorem 3.8 in \cite{Bona--02}.
\par

We have been tempted to state and prove 
a theorem about peripheral subgroups of $3$-manifolds
with boundary components of genus $\ge 2$.
But the situation is complicated,
because the boundary of an annulus in the characteristic surface $\mathcal A$
can be positioned in many ways relatively to $\partial M$.
Despite some effort, we did not arrive at a statement pleasing enough
to be stated here.
\par

To give some idea of the homotopy annulus theorem,
let us state it in the particular case of a manifold with boundary a union of tori,
a case for which no annuli are needed in the JSJ-like decomposition (see
\cite{NeSw--97}, Page 38).

\begin{Ahah}[\textbf{particular case}]
Consider the situation of Theorem \ref{JSJ},
and assume \emph{moreover} 
that $\partial M$ is a non-empty union of incompressible tori.
Assume that there exists an essential map (not necessariliy an embedding)
$\varphi : (\mathbf A, \partial \mathbf A) \longrightarrow (M, \partial M)$
of the annulus into $M$.
\par
Then there exists a homotopy
\begin{equation*}
\varphi_t \, : \,  (\mathbf A, \partial \mathbf A) \longrightarrow (M, \partial M) ,
\hskip.5cm 0 \le t \le 1 ,
\end{equation*}
such that $\varphi_0 = \varphi$ and
$\varphi_1$ is an essential map (not necessarily an embedding)
into a Seifert piece of the JSJ decomposition
(in particular, this decomposition has at least one Seifert piece).
\par
Moreover, in that Seifert piece, there exist embedded essential annuli
(vertical ones).
\end{Ahah}

\medskip

We cannot expect that $\varphi_1$ is homotopic to an embedding.
But, in a Seifert component, 
there are plenty of incompressible vertical annuli.
Hence we can strenghten the conclusion of the annulus theorem,
and conclude to the existence of an essential embedding
with vertical images inside some Seifert component.

\section{\textbf{Digression on the terminology of the literature}}
\label{sectionDigression}

\subsection{On essential annuli}
\label{On essential annuli}

Embedded annuli and tori play the key role in our arguments.
Since the literature concerning the related terminology 
is rather messy in our opinion,
we review the following definitions.
\par

Let $M$ be a bounded $3$-manifold and $A$ an annulus.
A mapping 
\begin{equation*}
f : (A, \partial A) \longrightarrow (M, \partial M)
\end{equation*}
is \textbf{W-essential}, or essential in the sense of Page 24 of \cite{Wald--78},
if the induced morphism of groups 
$\pi_1(A) \longrightarrow \pi_1(M)$
and the induced morphism of pointed sets 
$\pi_1(A, \partial A) \longrightarrow \pi_1(M, \partial M)$
are both injective ($f$ need not be an embedding).
Observe that $\pi_1(A, \partial A)$ has precisely two elements,
the base point and the non-trivial element represented by a spanning arc;
it follows that this definition of ``essential'' is equivalent to
that of \cite{Simo--76} (Page 206) or that of \cite{CaFe--76} (Page 220).
\par
A mapping $f : (A, \partial A) \longrightarrow (M, \partial M)$ is
\textbf{non-degenerate} 
if the homomorphism $\pi_1(A) \longrightarrow \pi_1(M)$ is injective
and if $f$ is not homotopic (as a map of pairs)
to a map $g : (A, \partial A) \longrightarrow (M, \partial M)$
with $g(A) \subset \partial M$
(we follow \cite{JaSh--79L}, Pages 121--122).
If $M$ is irreducible and if $\partial M$ is incompressible, 
Lemma IV.1.3 of \cite{JaSh--79L} shows that $f$ is non-degenerate in this sense
if and only if $f$ is $W$-essential.
\par
It seems that the terminology with ``essential'' becomes standard;
see for example \cite{Bona--02}, just before his Theorem 2.14.
\par
Thus, for an annulus properly embedded and incompressible 
in an irreducible manifold $M$ with boundary a union of tori,
we have a priori two notions: 
\begin{itemize}
\item[(i)]
it can be \textbf{$\partial$-incompressible}, 
or equivalently \textbf{not boundary parallel}
(see Proposition \ref{onannuli}), 
\item[(ii)]
it can be W-essential, or equivalently non-degenerate.
\end{itemize}
In fact, \textbf{these four notions are  equivalent}.
\par
Indeed, on the one hand, ``boundary parallel'' clearly implies ``degenerate''.
On the other hand, Lemma 5.3 of \cite{Wald--68} contains more than is necessary
to show that ``degenerate'' implies ``boundary parallel''.
Here is a weakened version of this Lemma 5.3, with Waldhausen's notation.

\begin{Lem}
Let $M$ be an irreducible $3$-manifold.
Let $G$ be an incompressible boundary component of $M$,
and let $F$ be an incompressible surface properly embedded in $M$
such that $\partial F \subset G$.
Suppose that there exists a homotopy
$H : F \times \mathbf I \longrightarrow M$
such that $H(F \times \{0\}) = F$
and $H(\partial(F \times \mathbf I) \smallsetminus (F \times \{0\}) \subset G$.
\par
Then $F$ is boundary parallel, and more precisely is parallel
to a surface contained in $G$.
\end{Lem}

\subsection{On the terminology of Raymond and Orlik}

For the reader who wishes to read \cite{Raym--68} and \cite{OrRa--68},
we offer the following dictionnary.
\begin{itemize}
\item[$\circ$]
$M$ is a compact and connected $3$-manifold on which $SO(2)$ acts.
For these authors, $M$ can be non-orientable
\emph{but we assume in this paper that $M$ is orientable},
so that their symbol $\epsilon$ takes always the value $o$
(small ``o'').
\item[$\circ$]
$M^*$ is the orbit space (our $\mathcal B$);
in our case, $M^*$ is a compact orientable surface,
and $g \ge 0$ is its genus.
\item[$\circ$]
$F$ is the fixed point set and $F^*$ is its homeomorphic image in $M^*$;
the number of connected components of $F$ 
is denoted by $h$ in \cite{Raym--68} and by $\hbar$ in \cite{OrRa--68}.
\item[$\circ$]
$E$ denotes the set of exceptional orbits;
its cardinal can be any non-negative integer.
$SE$ is the set of special exceptional orbits;
since $M$ is orientable here, $SE = \emptyset$,
so that its cardinal $t$ is always equal to $0$.
\end{itemize}
Hence, in our case, $M_{\epsilon, s, \hbar, t}$ is always
$M_{o, g, \hbar, 0}$.
\par
We do not need to comment on the Seifert invariants,
namely on $b$, which is a variant of the Euler class
(caution: there is a sign problem there),
and on $(\alpha_j, \beta_j)$, which are the usual Seifert invariants.
\par
The projective plane is denoted by $P$, the Klein bottle by $K$.
The ``non-orientable handle'' $N$,
which is the non-trivial $\mathbf S^2$-bundle over $\mathbf S^1$,
does not play any role for us.




\end{document}